\documentclass[graybox]{svmult}

\usepackage{amssymb,amsmath} 

\usepackage[polish,english]{babel}
\usepackage[T1]{fontenc}

\smartqed

\usepackage{mathptmx}       
\usepackage{helvet}         
\usepackage{courier}        
\usepackage{type1cm}        
\usepackage{makeidx}         
\usepackage{graphicx}        
\usepackage{multicol}        
\usepackage[bottom]{footmisc}

\makeindex             

\begin{document}

\title*{On some recent applications of stochastic convex ordering theorems to some functional inequalities for convex functions - a survey}
\titlerunning{Stochastic convex ordering theorems to  functional inequalities}
\author{Teresa Rajba}
\institute{Teresa Rajba \at University of Bielsko-Biala, Willowa 2, 43-309 Bielsko-Biala, Poland, \email{trajba@ath.bielsko.pl}
}
%
%

\maketitle

\abstract{This is a survey paper concerning some theorems on stochastic convex ordering and their applications to functional inequalities for convex functions. We present the recent results on those subjects.\smallskip}

\section{Introduction}
\label{introv2}
In the present paper we look at Hermite-Hadamard type inequalities from the perspective provided by the stochastic convex order.
 This approach is mainly due to Cal and C\'{a}rcamo. In the paper \cite{terCalCarc}, the Hermite-Hadamard type inequalities are interpreted in terms of the convex stochastic ordering between random variables. Recently, also in \cite {terFlorea, terOlbrysSzostok2014, terRajba14, terRajba12f, terRajba12g, terRajba2016a, terSzostok2014a, terTSzostok2014, terTSzostok2016b}), the Hermite-Hadamard inequalities are
studied  based on the convex ordering properties. Here we want to attract the readers attention
to some selected topics by presenting some theorems on the convex ordering that can be useful in the study of the Hermite-Hadamard type inequalities.

The Ohlin lemma \cite{terOhlin69} on sufficient conditions for convex stochastic ordering was 
 first used in \cite{terRajba12f},
to get a simple proof of some known Hermite-Hadamard type inequalities as well as to
obtaining new Hermite-Hadamard type inequalities. In \cite{terOlbrysSzostok2014, terTSzostok2014, terTSzostok2016b}, the authors used the Levin-Ste\v{c}kin theorem \cite{terLevinSteckin1960} to study Hermite-Hadamard type inequalities.

Many results on higher order generalizations of the Hermite-Hadamard type inequality one can found, among others, in \cite{terBes08, terBesPal02, terBesPal03, terBesPal04, terBesPal10, terDraPe2000, terRajba12f, terRajba12g}. In recent papers \cite{terRajba12f, terRajba12g} the theorem of M. Denuit, C.Lef\`{e}vre and  M. Shaked \cite{terDenLefSha98}  was used to prove Hermite-Hadamard type inequalities for higher-order convex functions. The theorem of M. Denuit, C.Lef\`{e}vre and  M. Shaked \cite{terDenLefSha98} on sufficient conditions for $s$-convex ordering is a counterpart of the Ohlin lemma concerning convex ordering. A theorem on necessary and sufficient conditions for higher order convex stochastic ordering, which is a counterpart of the Levin-Ste\v{c}kin theorem \cite{terLevinSteckin1960} concerning convex stochastic ordering, is given in the paper \cite{terRajba2016a}. Based on this theorem, useful criteria for the verification of higher order convex stochastic ordering are given. These criteria can be useful in the study of Hermite-Hadamard type inequalities for higher order convex functions, and in particular inequalities between the quadrature operators. They may be easier to verify the higher order convex orders, than those given in \cite{terDenLefSha98, terKARLNOVIK}.

In Section \ref{tergeneralization}, we give simple proofs of known as well as new Hermite-Hadamard type inequalities, using Ohlin's Lemma and the Levin-Ste\v{c}kin theorem.  

In Sections \ref{terdifferentiation} and \ref{terdifferentiation two}, we study inequalities of the Hermite-Hadamard type involving numerical differentiation formulas of the first order and the second order, respectively.

In Section \ref{tern order}, we give simple proofs of Hermite-Hadamard type inequalities for higher-order convex functions, using the theorem of M. Denuit, C.Lef\`{e}vre and  M. Shaked, and a generalization of the Levin-Ste\v{c}kin theorem to higher orders. These results are applied to derive some inequalities between quadrature operators.

\section{Some generalizations of the Hermite-Hadamard inequality}
\label{tergeneralization}
Let $f \colon [a,b]\to\mathbb{R}$ be a convex function ($a,b\in \mathbb{R},$ $a<b$). The following double inequality
\begin{equation}
f\left(\frac{a+b}2\right)\leq \frac1{b-a}\int_a^b f(x)\,dx \leq \frac{f(a)+f(b)}2 \label{tereq:oli0}
\end{equation}
is known as the Hermite-Hadamard inequality (see  \cite{terDraPe2000} for many generalizations and applications of \eqref{tereq:oli0}).

In many papers, the Hermite-Hadamard type inequalities are studied  based on the convex stochastic ordering properties (see, for example, \cite {terFlorea, terOlbrysSzostok2014, terRajba14, terRajba12f, terRajba12g, terSzostok2014a, terTSzostok2014}). In the paper \cite{terRajba12f}, the Ohlin lemma  on sufficient conditions for convex stochastic ordering  is used to get a simple proof of some known Hermite-Hadamard type inequalities as well as to obtain new Hermite-Hadamard type inequalities. Recently, the Ohlin lemma is also used  to study the inequalities of the Hermite-Hadamard type for convex functions in \cite {terOlbrysSzostok2014, terRajba14,  terSzostok2014a, terTSzostok2014}. In \cite {terRajba12g}, also the inequalities of the Hermite-Hadamard type for delta-convex functions are studied by using the Ohlin lemma.
In the papers \cite {terOlbrysSzostok2014, terSzostok2014a, terTSzostok2014}, furthermore, the Levin-Ste\v{c}kin theorem \cite {terLevinSteckin1960} (see also \cite {terNicPer06}) is used  to examine the Hermite-Hadamard type inequalities. This theorem gives necessary and sufficient conditions for the stochastic convex ordering.

Let us recall some basic notions and results on the stochastic convex order (see, for example, \cite {terDenLefSha98}). As usual, $F_X$ denotes the distribution function of a random variable $X$ and $\mu_X$ is the distribution corresponding to $X$. For real valued random variables $X,Y$ with a finite expectation, we say that $X$ is dominated by $Y$ in \textit{convex ordering} sense, if
 $$\mathbb{E} f(X) \leq \mathbb{E} f(Y)$$
for all convex functions $f \colon \mathbb{R} \to \mathbb{R}$ (for which the expectations exist). In that case we write $X \leq_{\textit{cx}} Y$, or $\mu_X \leq_{\textit{cx}} \mu_Y$.
\par\bigskip
In the following Ohlin's lemma \cite{terOhlin69}, are given sufficient conditions for convex stochastic ordering.
\begin{lemma}[Ohlin \cite{terOhlin69}]\label{terlemma:3}
Let $X,Y$ be two random variables such that $\mathbb{E} X=\mathbb{E} Y$. If the distribution functions $F_X, F_Y$ cross exactly one time, i.e., for some $x_0$ holds
\begin{equation*}
F_X(x) \leq F_Y(x) \text{ if } x < x_0 \quad\text{and}\quad  F_X(x) \geq F_Y(x) \text{ if } x > x_0,
\end{equation*}
then 
\begin{equation}\label{tereq:levin1}
\mathbb{E} f(X) \leq \mathbb{E} f(Y)
\end{equation}
for all convex functions $f \colon \mathbb{R} \to \mathbb{R}$.
\end{lemma}
\par\bigskip
The inequality \eqref{tereq:oli0} may be easily proved with the use of the Ohlin lemma (see\cite{terRajba12f}). Indeed, let $X$, $Y$, $Z$ be three random variables with the distributions $\mu_X = \delta_{(a+b)/2}$, $\mu_Y$ which is equally distributed in $[a,b]$ and $\mu_Z = \frac{1}{2}(\delta_a+\delta_b)$, respectively. Then it is easy to see that the pairs $(X,Y)$ and $(Y,Z)$ satisfy the assumptions of the Ohlin lemma, and using \eqref{tereq:levin1}, we obtain \eqref{tereq:oli0}.

\par\bigskip
Let $a<c<d<b$. Let $f \colon I \to \mathbb{R}$ be a convex function, $a,b \in I$. Then (see \cite{terHardy})
\begin{equation}
\frac{f(c)+f(d)}{2}-f\left(\frac{c+d}{2}\right)\leq \frac{f(a)+f(b)}{2} - f\left(\frac{a+b}{2}\right). \label{tereq:ID}
\end{equation}
\par\bigskip

To prove \eqref{tereq:ID} from the Ohlin lemma, it suffices to take random variables $X,Y$ (see \cite{terMihai}) with
 $$\mu_X= \frac{1}{4}\left(\delta_c+\delta_d\right)+\frac{1}{2}
 \delta_{(a+b)/2},$$
  $$\mu_Y= \frac{1}{4}\left(\delta_a+\delta_b\right)+\frac{1}{2}\delta_{\delta_{(c+d)/2}}.$$   
Then, by Lemma \ref{terlemma:3}, we obtain
\begin{equation}
\frac{f(c)+f(d)}{2}+f\left(\frac{a+b}{2}\right)\leq \frac{f(a)+f(b)}{2} +f\left(\frac{c+d}{2}\right),  \label{eq:IDa}
\end{equation}
which implies \eqref{tereq:ID}.

\par\bigskip
Similarly, it can be proved the Popoviciu inequality
\begin{equation}
\frac{2}{3}\left[ f\left(\frac{x+y}2\right)+f\left(\frac{y+z}2\right)+f\left(\frac{z+x}2\right)\right]\leq \frac{f(x)+f(y)+f(z)}3 + f\left(\frac{x+y+z}{3}\right),
\label{tereq:pop} 
\end{equation}
where $x,y,z \in I$ and $f \colon I \to \mathbb{R}$ is a convex function. To prove \eqref{tereq:pop} from the Ohlin lemma, it suffices (assuming $x\leq y \leq z$ ) to take random variables $X,Y$ (see \cite{terMihai}) with
$$
\mu_X= \frac{1}{4}\left(\delta_{(x+y)/2}+\delta_{(y+z)/2}+\delta_{(z+x)/2}\right),
$$
$$
\mu_Y= \frac{1}{6}\left(\delta_x+\delta_y+\delta_z\right)+\frac{1}{2}
 \delta_{(x+y+z)/3}.
$$

\par\bigskip
Convexity has a nice probabilistic characterization, known as Jensen's inequality (see \cite{terBillingsley1995}).
\begin{proposition}[\cite{terBillingsley1995}]
A function $f \colon (a,b) \to \mathbb{R}$ is convex if, and only if,
\begin{equation}
f(\mathbb{E}X)\leq  \mathbb{E}f(X)
\label{tereq:r1}
\end{equation}
for all $(a,b)$-valued integrable random variables $X$.
\end{proposition}
\par\bigskip

To prove \eqref{tereq:r1} from the Ohlin lemma, it suffices to take a random variable $Y$ (see \cite{terRajba14}) with
 $$\mu_Y= 
 \delta_{\mathbb{E}X},$$
then we have
\begin{equation}
\mathbb{E}f(Y)=f(\mathbb{E}X).
\label{tereq:r2}
\end{equation}

By the Ohlin lemma, we obtain $\mathbb{E}f(Y)\leq  \mathbb{E}f(X)$, then taking into account \eqref{tereq:r2}, this implies \eqref{tereq:r1}.
\par\bigskip
\begin{remark}
Note, that in \cite{terMroRajWas}, the Ohlin lemma was used to obtain a solution of the problem of Ra\c{s}a concerning inequalities for Bernstein operators.
\end{remark}
\par\bigskip

In \cite{terFe}, Fej\'er gave a generalization of the inequality \eqref{tereq:oli0}.

\begin{proposition}[\cite{terFe}]\label{terprop:1}
Let $f \colon I\to\mathbb{R}$ be a convex function defined on a real interval $I$, $a,b \in I$ with $a<b$ and let $g \colon [a,b]\to\mathbb{R}$ be non negative and symmetric with respect to the point $(a+b)/2$ (the existence of integrals is assumed in all formulas). Then
\begin{equation}
f\left(\frac{a+b}2\right)\cdot\int_a^b g(x)\,dx\leq\int_a^b f(x)g(x)\,dx\leq\frac{f(a)+f(b)}2\cdot\int_a^bg(x)\,dx. \label{tereq:1}
\end{equation}
\end{proposition}

The double inequality \eqref{tereq:1} is known in the literature as the Fej\'er inequality or the  Hermite-Hadamard-Fej\'er inequality (see \cite{terDraPe2000, terMiLa, terPecaricProschanTong1992} for the historical background).
\begin{remark}[\cite{terRajba12f}]
Using the Ohlin lemma (Lemma \ref{terlemma:3}), we get a simple proof of \eqref{tereq:1}. Let $f$ and $g$ satisfy the assumptions of Proposition \ref{terprop:1}. Let $X$, $Y$, $Z$ be three random variables such that $\mu_X = \delta_{(a+b)/2}$, $\mu_Y(dx) = ( \int_a^b g(x) dx )^{-1}g(x)dx$, $\mu_Z = \frac{1}{2}(\delta_a+\delta_b)$. Then, by Lemma \ref{terlemma:3}, we obtain that $X \leq_{\textit{cx}} Y$ and $Y \leq_{\textit{cx}} Z$, which implies \eqref{tereq:1}.
\end{remark}
\begin{remark}\label{terrem:2}

Note that for $g(x) = w(x)$ such that $\int_{a}^b w(x) dx = 1$, the inequality \eqref{tereq:1} can be rewritten in the form
	\begin{equation}
	f\left(\frac{a+b}{2}\right) \leq \int_a^b f(x) w(x) dx \leq \frac{f(a)+f(b)}{2}.
	\label{tereq:2}
	\end{equation}
	

Conversely, from the inequality \eqref{tereq:2}, it follows \eqref{tereq:1}. Indeed, if 
$\int_a^b g(x) dx>0$, it suffices to take $w(x) = \left(\int_a^b g(x) dx\right)^{-1} g(x)$. If $\int_a^b g(x) dx=0$, then \eqref{tereq:1} is obvious.
\end{remark}

For various modifications of \eqref{tereq:oli0} and \eqref{tereq:1} see e.g. \cite{terBesPal03, terBesPal04, terBesPal10, terCzinder06, terCzinderPal04, terDraPe2000}, and the references given there.

As Fink noted in \cite{terFink98}, one wonders what the symmetry has to do with the inequality \eqref{tereq:1} and if such an inequality holds for other functions (cf. \cite[p. 53]{terDraPe2000}).

As an immediate consequence of Lemma \ref{terlemma:3}, we obtain the following theorem, which is a generalization of the Fej\'er inequality.

\begin{theorem}[\cite{terRajba12f}]\label{terlemma:4}
Let $0 < p <1$. Let $f \colon I \to \mathbb{R}$ be a convex function, $a,b \in I$ with $a<b$. Let $\mu$ be a finite measure on $\mathcal{B}([a,b])$ such that (i) $\mu([a,pa+qb]) \leq pP_0$, (ii) $\mu((pa+	qb,b])\leq qP_0$, (iii) $\int_{[a,b]} x \mu(dx) = (pa+qb) P_0$, where $q=1-p$, $P_0 = \mu([a,b])$. Then
\begin{equation}
f(pa+qb)P_0 \leq \int_{[a,b]} f(x) \mu (dx) \leq [pf(a) + q f(b)] P_0.
\label{tereq:3}
\end{equation}
\end{theorem}

Fink proved in \cite{terFink98} a general weighted version of the Hermite-Hadamard inequality. In particular, we have the following probabilistic version of this inequality.

\begin{proposition}[\cite{terFink98}]\label{terprop:Fink}
Let $X$ be a random variable taking values in the interval $[a,b]$ such that $m$ is the expectation of $X$ and $\mu_X$ is the distribution corresponding to $X$. Then 
\begin{equation}
f\left(m\right)\leq\int_a^b f(x)\:\mu_X(dx)\leq\frac{b-m}{b-a} f(a)+ \frac{m-a}{b-a} f(b). \label{tereq:f}
\end{equation}
\end{proposition}


\par\bigskip
Moreover, in \cite{terFlorea} it was proved that, starting from such a fixed random variable $X$, we can fill the whole space
between the Hermite-Hadamard bounds by highlighting some parametric families of random variables. The authors propose two alternative constructions based on the convex ordering
properties.
\par\bigskip
In \cite{terRajba14}, based on Lemma \ref{terlemma:3}, a very simple proof of Proposition \ref{terprop:Fink} is given. Let $X$ be a random variable satisfying the assumptions of Proposition \ref{terprop:Fink}. Let $Y$, $Z$ be two random variables such that $\mu_Y = \delta_{m}$,   $\mu_Z = \frac{b-m}{b-a}\delta_a+\frac{m-a}{b-a}\delta_b)$. Then, by Lemma \ref{terlemma:3}, we obtain that $Y\leq_{\textit{cx}} X$ and $X \leq_{\textit{cx}} Z$, which implies \eqref{tereq:f}.

In \cite{terRajba12f}, some results related to the Brenner-Alzer inequality  are given. In the paper \cite{terBakPecPer12} by M. Klari\v{c}i\'{c} Bakula, J. Pe\v{c}ari\'{c} and J.  Peri\'{c}, some improvements of various forms of the Hermite-Hadamard inequality can be found; namely, that of Fej\'{e}r, Lupas, Brenner-Alzer, Beesack-Pe\v{c}ari\'{c}. These improvements imply the Hammer-Bullen inequality.
In 1991, Brenner and Alzer \cite{terBrenAlz91} obtained the following result generalizing Fej\'er's result as well as the result of Vasi\'c and Lackovi\'c (1976) \cite{terVasLac76} and Lupas (1976) \cite{terLu76} (see also \cite{terPecaricProschanTong1992}).

\begin{proposition}[\cite{terBrenAlz91}]\label{terprop:10}
Let $p,q$ be given positive numbers and $a_1 \leq a < b \leq b_1$. Then the inequalities
\begin{equation}
f\left(\frac{pa+qb}{p+q}\right) \leq \frac{1}{2y} \int_{A-y}^{A+y} f(t) dt \leq  \frac{pf(a)+qf(b)}{p+q}
\label{tereq:9}
\end{equation}
hold for $A=\frac{pa+qb}{p+q}$, $y>0$, and all continuous convex functions $f\colon [a_1,b_1] \to\mathbb{R}$ if, and only if,
$$
y \leq \frac{b-a}{p+q} \min\{p,q\}.
$$
\end{proposition}
\begin{remark}\label{terremark:11}
It is known \cite[p. 144]{terPecaricProschanTong1992} that  under the same conditions Hermite-Hadamard's inequality holds, the following refinement of \eqref{tereq:9}:
\begin{equation}
f\left(\frac{pa+qb}{p+q}\right) \leq \frac{1}{2y} \int_{A-y}^{A+y} f(t)dt \leq \frac{1}{2} \left\{ f(A-y) + f(A+y)\right\} \leq \frac{pf(a)+qf(b)}{p+q}
\label{tereq:10}
\end{equation}
holds.
\end{remark}

In the following theorem we give some generalization of the Brenner and Alzer inequalities \eqref{tereq:10}, which we prove using the Ohlin lemma.
\begin{theorem}[\cite{terRajba12f} ]\label{terth:12}
Let $p,q$ be given positive numbers, $a_1 \leq a < b \leq b_1$, $0 < y \leq \frac{b-a}{p+q}\min\{p,q\}$ and let $f\colon [a_1,b_1]\to\mathbb{R}$ be a convex function. Then
$$
f\left(\frac{pa+qb}{p+q}\right)\leq 
$$ 
$$\frac{\alpha}{2} \left\{ f(A-(1-\alpha)y) + f(A+(1-\alpha)y)\right\} + \frac{1}{2y} \int_{A-(1-\alpha)y}^{A+(1-\alpha)y} f(t)dt \leq 
$$
$$\frac{\alpha}{2n} \sum_{k=1}^n \left\{ f\left(A-y+k\frac{\alpha y}{n}\right) + f\left(A+y-k\frac{\alpha y}{n}\right)\right\} + \frac{1}{2y} \int_{A-(1-\alpha)y}^{A+(1-\alpha)y} f(t)dt\leq 
$$
\begin{equation}\frac{1}{2y} \int_{A-y}^{A+y} f(t) dt, \label{tereq:11}
\end{equation}
where $0 \leq \alpha \leq 1$, $n=1,2,\ldots$,
\begin{eqnarray}
\frac{1}{2y} \int_{A-y}^{A+y} f(t)dt & \leq & \frac{\beta}{2}\{f(A-y) + f(A+y)\} + (1-\beta)\frac{1}{2y} \int_{A-y}^{A+y} f(t)dt \nonumber \\
 & \leq & \frac{1}{2}\{ f(A-y)+f(A+y)\}, \label{tereq:12}
\end{eqnarray}
where $0 \leq \beta \leq 1$,
$$\frac{1}{2}\{ f(A-y) + f(A+y) \} \leq 
$$ 
$$(\frac{1}{2}-\gamma)\{ f(A-y-c) + f(A+y+c) \} + \gamma\{ f(A-y) + f(A+y) \} \leq
$$
\begin{equation} \frac{pf(a)+qf(b)}{p+q}, \label{tereq:13}\end{equation}
where $c = \min\{ b-(A+y), (A-y)-a \}$, $\gamma = \left|\frac{1}{2}-p\right|$.
\end{theorem}
\par\bigskip

To prove this theorem, it suffices to consider random variables $X$, $Y$, $W$, $Z$, $\xi_n$, $\eta$ and $\lambda$ such that:
\begin{eqnarray*}
\mu_X & = & \delta_{\frac{ pa + qb}{p+q}}, \\
\mu_Y(dx) & = & \frac{1}{2y} \chi_{[A-y,A+y]}(x)dx, \\
\mu_Z & = & \frac{p}{p+q}\delta_a + \frac{q}{p+q}\delta_b, \mu_W = \frac{1}{2}\delta_{A-y}+\frac{1}{2}\delta_{A+y}, \\
\mu_{\xi_n}(dx) & = & \frac{\alpha}{2n}\sum_{k=1}^n \{\delta_{A-y + k\frac{\alpha y}{n}} + \delta_{A+y-k\frac{\alpha y}{n}}\} + \frac{1}{2y} \chi_{[A-(1-\alpha)y, A+(1-\alpha)y]}(x)dx, \\
\mu_\eta(dx) & = & \frac{\beta}{2} \{ \delta_{A-y} + \delta_{A+y} \} + \frac{1-\beta}{2y}\chi_{[A-y,A+y]}(x)dx,\\
\mu_\lambda & = & (\frac{1}{2} - \gamma)\{ \delta_{A-y-c} + \delta_{A+y+c}\} + \gamma\{\delta_{A-y}+\delta_{A+y}\}.
\end{eqnarray*}
Then, using the Ohlin lemma, we obtain:
\begin{itemize}
	\item $X \leq_{\textit{cx}} Y$, $Y \leq_{\textit{cx}} W$ and $W \leq_{\textit{cx}} Z$, which implies the inequalities \eqref{tereq:10},
	\item $X \leq_{\textit{cx}} \xi_1$, $\xi_1 \leq_{\textit{cx}} \xi_n$ and $\xi_n \leq_{\textit{cx}} Y$, which implies \eqref{tereq:11},
	\item $Y \leq_{\textit{cx}} \eta$ and $\eta \leq_{\textit{cx}} W$, which implies \eqref{tereq:12},
	\item $W \leq_{\textit{cx}} \lambda$ and $\lambda \leq_{\textit{cx}} Z$, which implies \eqref{tereq:13}.
\end{itemize}
\begin{theorem}[\cite{terRajba12f} ]\label{terth:13}
Let $p,\:q$ be given positive numbers, $0 < \alpha < 1$, $a_1 \leq a < b \leq b_1$, $0< y \leq \frac{b-a}{p+q} \min \{ p, q\}$ and $0 \leq \frac{\alpha}{1-\alpha}y \leq \frac{b-a}{p+q} \min\{p,q\}$. Let $f \colon [a_1,b_1]\to\mathbb{R}$ be a~convex function. Then
\begin{eqnarray}
f(A) & \leq & \frac{\alpha}{y} \int_{A-y}^{A} f(t)dt+ \frac{(1-\alpha)^2}{\alpha y} \int_{A}^{A+\frac{\alpha}{1-\alpha}y} f(t)dt \nonumber \\
 & \leq & \alpha f(A-y) + (1-\alpha)f(A+\frac{\alpha}{1-\alpha}y) \nonumber \\
 & \leq & \frac{p}{p+q} f(a) + \frac{q}{p+q} f(b), \label{tereq:14}
\end{eqnarray}
where $A = \frac{pa +qb}{p+q}$.
\end{theorem}

Let $X$, $Y$, $Z$ and $W$ be random variables such that:
\begin{eqnarray*}
\mu_X & = & \delta_A ,\\
\mu_Y(dx) & = & \frac{\alpha}{y}\chi_{[A-y,A]}(x)dx + \frac{(1-\alpha)^2}{\alpha y} \chi_{[A,A+\frac{\alpha}{1-\alpha}y]}(x)dx, \\
\mu_W & = & \alpha\delta_{A-y} + (1-\alpha)\delta_{A+\frac{\alpha}{1-\alpha}y}, \\
\mu_Z & = & \frac{p}{p+q}\delta_a + \frac{q}{p+q}\delta_b.
\end{eqnarray*}
Then, using the Ohlin lemma, we obtain $X \leq_{\textit{cx}} Y$, $Y \leq_{\textit{cx}} W$, $W \leq_{\textit{cx}} Z$, which implies the inequalities \eqref{tereq:14}.
\begin{remark}\label{terremark:14}
If we choose $\alpha = \frac{1}{2}$ in Theorem \ref{terth:13}, then the inequalities \eqref{tereq:14} reduce to the inequalities \eqref{tereq:12}.
\end{remark}

\begin{remark}\label{terremark:15}
If we choose $\alpha = \frac{p}{p+q}$ and $y=(1-p)z$ in Theorem \ref{terth:13}, then we have
\begin{eqnarray*}
f(A) & \leq & \frac{p}{qz} \int_{A-\frac{q}{p+q}z}^A f(t)dt + \frac{q}{pz}\int_{A}^{A+\frac{p}{p+q}z} f(t)dt \\
 & \leq & \frac{p}{p+q} f(A-\frac{q}{p+q}z) + \frac{q}{p+q} f(A + \frac{p}{p+q}z) \\
 & \leq & \frac{p}{p+q} f(a) + \frac{q}{p+q}f(b),
\end{eqnarray*}
where $A = \frac{pa+qb}{p+q}$, $0 < z \leq b-a$.
\end{remark}

In the paper \cite{terSzostok2014a}, the author used Ohlin's lemma to prove some new inequalities of the Hermite-Hadamard type, which are a generalization of known Hermite-Hadamard type inequalities.

\begin{theorem}[\cite{terSzostok2014a}]
\label{ter1}
The inequality 
\begin{equation}
\label{ter2hh}
af(\alpha x+(1-\alpha )y)+(1-a)f(\beta x+(1-\beta)y)\leq \frac{1}{y-x}\int_{x}^yf(t)dt,
\end{equation}
with some $a,\alpha,\beta\in[0,1],$ $\alpha>\beta$
is satisfied for all $x,y\in\mathbb{R}$ and all continuous and convex functions $f:[x,y]\to\mathbb{R}$ if, and only if, 
\begin{equation}
\label{ter12}
a\alpha+(1-a)\beta=\frac12,
\end{equation}
and one of the following conditions holds true:

(i) $a+\alpha\leq 1$, 

(ii)  $a+\beta\geq 1$, 

(iii) $a+\alpha>1,$ $a+\beta<1$ and $a+2\alpha\leq 2.$
\end{theorem}

\begin{theorem}[\cite{terSzostok2014a}]
\label{terr3}
 Let  $a,b,c,\alpha\in(0,1)$ be numbers such that $a+b+c=1$. Then the inequality
\begin{equation}
\label{ter3HH}
af(x)+bf(\alpha x+(1-\alpha)y )+cf(y)\geq \frac{1}{y-x}\int_{x}^yf(t)dt
\end{equation}
is satisfied for all $x,y\in\mathbb{R}$ and all continuous and convex functions $f:[x,y]\to\mathbb{R}$ if, and only if, 
\begin{equation}
\label{ter121}
b(1-\alpha)+c=\frac12
\end{equation}
and one of the following conditions holds true:

(i) $a+\alpha\geq1,$ 

(ii) $a+b+\alpha\leq 1,$ 

(iii) $a+\alpha<1,$ $a+b+\alpha>1$ and $2a+\alpha\geq 1.$
\end{theorem}

Note that the original Hermite-Hadamard inequality consists of two parts. We treated these cases
separately. However, it is possible to formulate a result containing both inequalities. 

\begin{corollary}[\cite{terSzostok2014a}]
If $a,\alpha,\beta\in(0,1)$ satisfy (\ref{ter12}) and one of the conditions $(i)$, $(ii)$, $(iii)$ 
of Theorem \ref{ter1}, then the inequality
$$af(\alpha x+(1-\alpha y)+(1-a)f(\beta x+(1-\beta)y)\leq \frac{1}{y-x}\int_{x}^yf(t)dt\leq$$
$$ (1-\alpha)f(x)+(\alpha-\beta)f(ax+(1-a)y)+\beta f(y)$$
is satisfied for all $x,y\in\mathbb{R}$ and for all continuous and convex functions $f:\mathbb{R}\to\mathbb{R}.$
\end{corollary}

As we can see, the Ohlin lemma is very useful, however, it is worth noticing that in the case of some inequalities, the distribution functions cross more than once. Therefore a simple application of the Ohlin lemma is impossible. 

In the papers \cite{terOlbrysSzostok2014, terTSzostok2014}, the authors used the Levin-Ste\v{c}kin theorem  \cite{terLevinSteckin1960} (see also \cite{terNicPer06}, Theorem 4.2.7), which gives necessary and sufficient conditions for convex ordering of functions with bounded variation, which are distribution functions of signed measures.

\begin{theorem}[Levin, Ste\v{c}kin \cite{terLevinSteckin1960}]\label{terth:1.3}
Let $a,b\in\mathbb{R}$, $a<b$ and let $F_1,F_2 \colon [a,b]\to\mathbb{R}$ be functions with bounded variation such that $F_1(a)=F_2(a)$. Then, 
in order that
\begin{equation}
\int_a^bf(x)dF_1(x)\leq\int_a^bf(x)dF_2(x)
\end{equation}
for all continuous convex functions $f \colon [a,b]\to\mathbb{R},$
it is necessary and sufficient that $F_1$ and $F_2$ verify the following three conditions:
\begin{eqnarray}
F_1(b)& = & F_2(b)\label{tereq:ls1.4},\\
\int_a^b F_1(x)dx & = & \int_a^b F_2(x)dx, \label{eq:ls1.6}\\
  \int_a^x F_1(t)dt & \leq & \int_a^x F_2(t)dt\quad \text{for all}\quad x\in(a,b).\label{tereq:ls1.5} 
\end{eqnarray}
\end{theorem}

Define the number of sign changes of a function $\varphi \colon \mathbb{R} \to \mathbb{R}$ by
$$
S^-(\varphi) = \sup \{S^-[\varphi(x_1),\varphi(x_2),\ldots,\varphi(x_k)]\colon x_1<x_2<\ldots x_k\in\mathbb{R},\: k\in\mathbb{N}\},
$$
where $S^-[y_1,y_2,\ldots,y_k]$ denotes the number of sign changes in the sequence $y_1$, $y_2$,$\ldots,$ $y_k$ (zero terms are being discarded). Two real functions $\varphi_1,\varphi_2$ are said to have $n$ crossing points (or cross each other $n$-times) if $S^-(\varphi_1-\varphi_2)=n$. Let $a=x_0<x_1< \ldots<x_n<x_{n+1}=b$. We say that the functions  $\varphi_1,\varphi_2$ crosses $n$-times at the points $x_1,x_2,\ldots,, x_n$ (or that $x_1,x_2,\ldots,, x_n$ are the points of sign changes of $\varphi_1-\varphi_2$) if $S^-(\varphi_1-\varphi_2)=n$ and there exist $a<\xi_1<x_1< \ldots<\xi_n<x_n<\xi_{n+1}<b$ such that $S^-[\xi_1,\xi_2,\ldots,\xi_{n+1}]=n$.

Szostok \cite{terTSzostok2014} used Theorem \ref{terth:1.3} to make an observation, which is more general than Ohlin's lemma and concerns the situation when the functions $F_1$ and $F_2$ have more crossing points than one. 
In \cite{terTSzostok2014} is given some useful modification of the Levin-Ste\v{c}kin theorem  \cite{terLevinSteckin1960}, which can be rewritten in the following form. 

\begin{lemma}[\cite{terTSzostok2014}]\label{terlemma:1.5}
Let $a,b\in\mathbb{R}$, $a<b$ and let $F_1,F_2 \colon (a,b)\to\mathbb{R}$ be functions with bounded variation such that $F(a)=F(b)=0$, $\int_a^b F(x)dx=0$, where $F=F_2-F_1$. Let $a<x_1< \ldots<x_m<b$ be the points of sign changes of the function $F$. Assume that $F(t)\geq 0$   for $t\in(a,x_1)$.  
\begin{itemize}
\item If $m$ is even then 
the inequality 
\begin{equation}\label{tereq:ord1}
\int_a^bf(x)dF_1(x)\leq\int_a^bf(x)dF_2(x)
\end{equation}
is not satisfied by all continuous convex functions $f \colon [a,b]\to\mathbb{R}.$ 

\item If $m$ is odd, define $A_i$ ($i=0,1,\ldots,m$, $x_0=a$, $x_{m+1}=b$)
$$
A_i = \int_{x_i}^{x_{}i+1} |F(x)|dx.
$$
Then the inequality \eqref{tereq:ord1} is satisfied for all continuous convex functions $f \colon [a,b]\to\mathbb{R},$  if, and only if, the following inequalities hold true:
\begin{equation}\label{tereq:a1}
\begin{split}
A_0 & \geq  A_1, \\
A_0+A_2 & \geq  A_1+A_3 ,\\
& \vdots \\
A_0+A_2+ \ldots+A_{m-3} & \geq  A_1+A_3+ \ldots+A_{m-2}.
\end{split}
\end{equation}
\end{itemize}
\end{lemma}

\begin{remark}[\cite{terRajba2016a}]\label{terrmk:a2}
Let
\begin{equation*}
H(x)=\int_a^x F(t)dt.
\end{equation*}
Then the inequalities  \eqref{tereq:a1}
are equivalent to the following inequalities
\begin{equation*}
H(x_2)\geq  0,\: H(x_4)\geq  0,\:H(x_6)\geq  0, \:\ldots, \:H(x_{m-1})\geq 0.
\end{equation*}
\end{remark}

In \cite{terTSzostok2014}, Lemma \ref{terlemma:1.5} is used to prove results, which extend the inequalities \eqref{ter2hh} and \eqref{ter3HH} and inequalities between quadrature operators. 
\begin{theorem}[\cite{terTSzostok2014}]
\label{terthlH}
 Let numbers $a_1,a_2,a_3,\alpha_1,\alpha_2,\alpha_3\in(0,1)$ satisfy
$a_1+a_2+a_3=1$ and $\alpha_1>\alpha_2>\alpha_3.$

Then the inequality 
\begin{equation}
\label{tertlH}
\sum_{i=1}^3a_if(\alpha_ix+(1-\alpha_i)y)\leq \frac{1}{y-x}\int_{x}^yf(t)
\end{equation}
is satisfied by all convex functions $f:[x,y]\to\mathbb{R}$ if, and only if, we have 
\begin{equation}
\label{terE1lev}
\sum_{i=1}^3a_i(1-\alpha_i)=\frac12
\end{equation}

 and one of the following conditions is satisfied

(i) $a_1\leq 1-\alpha_1$ and $a_1+a_2\geq 1-\alpha_3,$

(ii) $a_1\geq 1-\alpha_2$ and $a_1+a_2\geq 1-\alpha_3,$

(iii) $a_1\leq 1-\alpha_1$ and $a_1+a_2\leq 1-\alpha_2,$

(iv) $a_1\leq 1-\alpha_1, a_1+a_2\in(1-\alpha_2,1-\alpha_3)$ and $2\alpha_3\geq a_3,$

(v) $a_1\geq 1-\alpha_2,a_1+a_2<1-\alpha_3$ and $2\alpha_3\geq a_3,$

(vi) $a_1> 1-\alpha_1,a_1+a_2\leq 1-\alpha_2$ and $1-\alpha_1\geq\frac{a_1}{2},$

(vii) $a_1\in(1-\alpha_1,1-\alpha_2),$ $a_1+a_2\geq 1-\alpha_3,$ 
and $1-\alpha_1\geq\frac{a_1}{2}$,
(viii) $a_1\in(1-\alpha_1,1-\alpha_2),$ $a_1+a_2\in(1-\alpha_2,1-\alpha_3),$ 
 $1-\alpha_1\geq\frac{a_1}{2}$ and $2a_1(1-\alpha_1)+2a_2(1-\alpha_2)\geq (a_1+a_2)^2.$
\end{theorem}

To prove Theorem \ref{terthlH}, we note that, if the inequality \eqref{tertlH} is satisfied for every convex function $f$ defined on the interval $[0,1]$,
then it is satisfied by every convex function $f$ defined on a given interval $[x,y].$ Therefore, without loss of generality, it suffices to consider the interval $[0,1]$ in place of $[x,y].$

To prove Theorem \ref{terthlH}, we consider the functions $F_1,F_2:\mathbb R\to\mathbb{R}$ given by the following formulas
\begin{equation}
\label{terF11}
F_1(t):=\left\{\begin{array}{ll}
0,&\ t<1-\alpha_1,\\
a_1,&\ t\in[1-\alpha_1,1-\alpha_2),\\
a_1+a_2,&\ t\in[1-\alpha_2,1-\alpha_3),\\
1,&\ t\geq 1-\alpha_3,
\end{array}\right.
\end{equation}
and
\begin{equation}
\label{terF21}
F_2(t):=\left\{\begin{array}{ll}
0,&\quad t<0,\\
t,&\quad t\in[0,1),\\
1,&\quad t\geq 1.
\end{array}\right.
\end{equation}
Observe that the equality \eqref{terE1lev} gives us 
$$\int_{0}^1tdF_1(t)=\int_{0}^1tdF_2(t).$$
Further, it is easy to see that in the cases $(i),(ii)$ and $(iii)$ the pair $(F_1,F_2)$ 
crosses exactly once and, consequently, the inequality \eqref{tertlH} follows from the Ohlin lemma.

In the case $(iv)$, the pair $(F_1,F_2)$ crosses three times. Let $A_0,\dots,A_3$ be defined
as in Lemma \ref{terlemma:1.5}. In order to prove the inequality \eqref{tertlH}, we note that $A_0\geq A_1.$ However, since $A_0-A_1+A_2-A_3=0,$ we shall show that 
$A_2\leq A_3.$  We have 
$$A_2=\int_{a_1+a_2}^{1-\alpha_3}(t-a_1-a_2)dt=\frac{(1-\alpha_3-a_1-a_2)^2}{2}=\frac{a_3^2-2a_3\alpha_3+\alpha_3^2}{2}$$
and
$$A_3=\int_{1-\alpha_3}^1 (1-t)dt=\frac{\alpha_3^2}{2}.$$
This means that $A_2\leq A_3$ is equivalent to $2\alpha_3\geq a_3,$ as claimed.

We omit similar proofs in the cases $(v),(vi)$ and $(vii)$ and we pass to the case $(vii).$
In this case, the pair $(F_1,F_2)$ crosses five times. 
We have
$$A_0=\int_0^{1-\alpha_1}tdt=\frac{(1-\alpha_1)^2}{2}$$
and
$$A_1=\int_{1-\alpha_1}^{a_1}(a_1-t)dt=a_1(a_1-(1-\alpha_1))-\frac{a_1^2-(1-\alpha_1)^2}{2}=
\frac{[a_1-(1-\alpha_1)]^2}{2}.$$
This means, that the inequality $A_0\geq A_1$ is satisfied if, and only if, $1-\alpha_1\geq\frac{a_1}{2}.$

Further,
$$A_2=\int_{a_1}^{1-\alpha_2}(t-a_1)dt=\frac{(1-\alpha_2)^2-a_1^2}{2}-a_1(1-\alpha_2-a_1)$$
and 
$$A_3=\int_{1-\alpha_2}^{a_1+a_2}(a_1+a_2-t)dt=(a_1+a_2)(a_1+a_2-(1-\alpha_2))-\frac{(a_1+a_2)^2-(1-\alpha_2)^2}{2},$$
 therefore, the inequality $A_0+A_2\geq A_3+A_1$ is satisfied if, and only if,
$$(1-\alpha_1)^2+(1-\alpha_2-a_1)^2\geq (a_1-1-\alpha_1)^2+(a_1+a_2-1+\alpha_2)^2,$$
which, after some calculations, gives us the last inequality from $(vii).$
\par\bigskip

Using  assertions (i) and (vii) of Theorem \ref{terthlH}, it is easy to get the following example. 

\begin{example}[\cite{terTSzostok2014}]
\label{terpr1}
Let $x,y\in\mathbb{R},$ $\alpha\in(\frac12,1)$  and $a,b\in(0,1)$ be such that $2a+b=1.$ 
Then the inequality 
\begin{equation}
\label{tersymlH}
af(\alpha x+(1-\alpha)y)+bf\left(\frac{x+y}{2}\right)+af((1-\alpha)x+\alpha y)\leq\frac{1}{y-x}\int_x^yf(t)dt
\end{equation}
is satisfied by all convex functions $f:[x,y]\to\mathbb{R}$ if, and only if,  $a\leq 2-2\alpha.$
\end{example}

In the next theorem, we obtain inequalities, which extend the second of the Hermite-Hadamard inequalities.
\begin{theorem}[\cite{terTSzostok2014}]
\label{terthrH}
 Let numbers $a_1,a_2,a_3,a_4\in(0,1),\alpha_1,\alpha_2,\alpha_3,\alpha_4\in[0,1]$ satisfy
$a_1+a_2+a_3+a_4=1$ and $1=\alpha_1>\alpha_2>\alpha_3>\alpha_4=0.$

Then the inequality 
\begin{equation}
\label{terrH}
\sum_{i=1}^4a_if(\alpha_ix+(1-\alpha_i)y)\geq \frac{1}{y-x}\int_{x}^yf(t)
\end{equation}
is satisfied by all convex functions $f:[x,y]\to\mathbb{R}$ if, and only if, we have 
\begin{equation}
\label{terE4}
\sum_{i=1}^4a_i(1-\alpha_i)=\frac12
\end{equation}
and one of the following conditions is satisfied:

(i) $a_1\geq1-\alpha_2$ and $a_1+a_2\geq1-\alpha_3,$

(ii) $a_1+a_2\leq 1-\alpha_2$ and $a_1+a_2+a_3\leq 1-\alpha_3,$ 

(iii) $1-\alpha_2\leq a_1$ and $1-\alpha_3\geq a_1+a_2+a_3,$

(iv) $1-\alpha_2\leq a_1,$ $1-\alpha_3\in(a_1+a_2,a_1+a_2+a_3)$ and $\alpha_3\leq 2a_4,$ 

(v) $1-\alpha_2\geq a_1+a_2,a_1+a_2+a_3>1-\alpha_3$ and $\alpha_3\leq 2a_4,$

(vi) $a_1<1-\alpha_2,$ $a_1+a_2\geq 1-\alpha_3$ and $2a_1+\alpha_2\geq 1,$

(vii) $a_1<1-\alpha_2,a_1+a_2>1-\alpha_2,a_1+a_2+a_3\leq 1-\alpha_3$ and $2a_1+\alpha_2\geq 1,$

(viii) 
$1-\alpha_2\in(a_1,a_1+a_2),1-\alpha_3\in(a_1+a_2,a_1+a_2+a_3),2a_1+\alpha_2
\geq 1$ and $2a_1(1-\alpha_3)+2a_2(\alpha_2-\alpha_3)\geq(1-\alpha_3)^2.$
\end{theorem}

To prove Theorem \ref{terthrH}, we assume that $F_1:\mathbb R\to\mathbb{R}$ is the function given by the following formula
\begin{equation}
\label{refF12}
F_1(t):=\left\{\begin{array}{ll}
0,&\ t<0,\\
a_1,&\ t\in[0,1-\alpha_1),\\
a_1+a_2,&\ t\in[1-\alpha_1,1-\alpha_2),\\
a_1+a_2+a_3,&\ t\in[1-\alpha_2,1),\\
1,&\ t\geq 1.
\end{array}\right.
\end{equation}
and let $F_2$ be the function given by \eqref{terF21}.
In view of \eqref{terE4}, we have 
$$\int_{0}^1F_1(t)dt=\int_{0}^1F_2(t)dt.$$
In cases $(i)-(iii)$ there is only one crossing point of $(F_2,F_1)$  
and our assertion is a consequence of the Ohlin lemma.

In the cases $(iv)-(vii)$, the pair $(F_2,F_1)$ crosses three times and, therefore, we have to use Lemma \ref{terlemma:1.5}.

In the case $(iv)$, the inequality \eqref{terrH} is satisfied by all convex functions $f$ if, and only if, $A_0\geq A_1.$ Further, we know that 
$$A_0-A_1+A_2-A_3=0,$$
which implies that the inequality $A_0\geq A_1$ is equivalent to $A_3\geq A_2.$ Clearly, we have
\begin{equation}
\label{terA2}
\begin{split}
A_2=&\int_{1-\alpha_3}^{1-a_4}(F_1(t)-F_2(t))dt=(\alpha_3-a_4)(1-a_4)-\frac{(1-a_4)^2-(1-\alpha_3)^2}{2}\\
=&(\alpha_3-a_4)\left(1-a_4+\frac{2-(\alpha_3+a_4)}{2}\right)
\end{split}
\end{equation}
and 
\begin{equation}
\label{terA3}
A_3=\int_{1-a_4}^1 (t-(1-a_4))dt=\frac{1-(1-a_4)^2}{2}-(1-a_4)a_4
\end{equation}
 i.e. $A_3\geq A_2$ is equivalent to $\alpha_3\leq 2a_4.$

We omit similar reasoning in the cases $(v),(vi)$ and $(vii)$ and we pass to the most interesting case $(viii).$
In this case, $(F_2,F_1)$ has 5 crossing points and, therefore, we must check that the inequalities 
$$A_0\geq A_1\;\;\textrm{and}\;\;A_0-A_1+A_2\geq A_3$$ are equivalent to
the inequalities of the condition $(viii)$, respectively. 
To this end, we write
$$A_0=\int_{0}^{a_1}(a_1-t)dt=\frac{a_1^2}{2},$$
$$A_1=\int_{a_1}^{1-\alpha_1}(t-(a_1+a_2))dt=\frac{(a_1+a_2-1+\alpha_1)^2}{2},$$
which means that $A_0\geq A_1$ if, and only if, $2a_1+\alpha_2\geq 1.$ 
Further, $A_2$ and $A_3$ are given by formulas \eqref{terA2} and \eqref{terA3}. Thus,
$A_0-A_1+A_2\geq A_3$ is equivalent to 
$$a_1^2+(a_1+a_2-(1-\alpha_2))^2\geq(1-\alpha_2-a_1)^2+(1-\alpha_3-a_1-a_2)^2,$$
which yields
$$2a_1(1-\alpha_3)+2a_2(\alpha_2-\alpha_3)\geq(1-\alpha_3)^2.$$
\par\bigskip

Using assertions $(ii)$ and $(vii)$ of Theorem \ref{terthrH}, 
we get the following example.
\begin{example}[\cite{terTSzostok2014}]
Let $x,y\in\mathbb{R},$ let $\alpha\in(\frac12,1)$ and let $a,b\in(0,1)$ be such that $2a+2b=1.$ 
Then, the inequality 
$$af(x)+bf(\alpha x+(1-\alpha)y)+bf\left((1-\alpha)x+\alpha y\right)+af(y)\geq\frac{1}{y-x}\int_x^yf(t)dt$$
is satisfied by all convex functions $f:[x,y]\to\mathbb{R}$ if, and only if, 
$a\geq\frac{1-\alpha}{2}.$
\end{example}

In the next theorem we show, that the same tools may be used to obtain some inequalities between quadrature operators, which do not involve the integral mean. 
\begin{theorem}[\cite{terTSzostok2014}]
\label{terthqo}
 Let $a,\alpha_1,\alpha_2,\beta\in(0,1)$ and let $b_1,b_2,b_3\in(0,1)$ satisfy
$b_1+b_2+b_3=1.$

Then, the inequality 
\begin{equation*}
af(\alpha_1x+(1-\alpha_1)y)+(1-a)f(\alpha_2x+(1-\alpha_2)y)\leq 
\end{equation*}
\begin{equation}
\label{terinqo}
 b_1f(x)+b_2f(\beta x+(1-\beta)y)+b_3f(y)
\end{equation}
is satisfied by all convex functions $f:[x,y]\to\mathbb{R}$ if, and only if, we have 
\begin{equation}
\label{terE}
b_2(1-\beta)+b_3=a(1-\alpha_1)+(1-a)(1-\alpha_2)
\end{equation}

 and one of the following conditions is satisfied:

(i) $a\leq b_1,$

(ii) $a\geq b_1+b_2,$

(iii) $\alpha_2\geq\beta$

or

(iv) $a\in(b_1,b_1+b_2),$ $\alpha_2<\beta$ and $(1-\alpha_1)b_1\geq(\alpha_1-\beta)(a-b_1).$

\end{theorem}
 
Now, using this theorem, we shall present positive and negative examples of inequalities of the type  \eqref{terinqo}.
\begin{example}[\cite{terTSzostok2014}] 
Let $\alpha\in\left(\frac12,1\right).$
The inequality 
$$\frac{f(\alpha x+(1-\alpha)y)+f((1-\alpha)x+\alpha y)}{2}\leq\frac{f(x)+f\left(\frac{x+y}{2}\right)+f(y)}{3}$$
is satisfied by all convex functions $f:[x,y]\to\mathbb{R}$ if, and only if, $\alpha\leq\frac56.$
\end{example}

\begin{example}[\cite{terTSzostok2014}]
 Let $\alpha\in\left(\frac12,1\right).$
The inequality 
$$\frac{f(\alpha x+(1-\alpha)y)+f((1-\alpha)x+\alpha y)}{2}\leq\frac16f(x)+\frac{2}{3}f\left(\frac{x+y}{2}\right)+\frac16f(y)$$
is satisfied by all convex functions $f:[x,y]\to\mathbb{R}$ if, and only if, $\alpha\leq\frac23.$
\end{example}

\section{Inequalities of the Hermite-Hadamard type involving numerical differentiation formulas of the first order}
\label{terdifferentiation}
In the paper \cite{terOlbrysSzostok2014}, expressions connected with numerical differentiation formulas of order 
$1$ are studied. The authors used the Ohlin lemma and the Levin-Ste\v{c}kin
theorem to study inequalities of the Hermite-Hadamard type connected with these expressions.

First, we recall the classical Hermite-Hadamard inequality    
\begin{equation}
\label{terHH}
f\left(\frac{x+y}{2}\right)\leq\frac{1}{y-x}\int_{x}^yf(t)dt\leq\frac{f(x)+f(y)}{2}.
\end{equation}
Now, let us write \eqref{terHH} in the form    
\begin{equation}
\label{terHH1}
f\left(\frac{x+y}{2}\right)\leq\frac{F(y)-F(x)}{y-x}\leq\frac{f(x)+f(y)}{2}.
\end{equation}
Clearly, this inequality is satisfied by every convex function $f$ and its primitive function $F$.
However,  \eqref{terHH1} may be viewed as an inequality involving two types of expressions used, in numerical integration and differentiation, respectively. Namely, 
$f\left(\frac{x+y}{2}\right)$ and $\frac{f(x)+f(y)}{2}$ are the simplest quadrature formulas 
used to approximate the definite integral, whereas 
$\frac{F(y)-F(x)}{y-x}$ is the simplest expression used to approximate the derivative of $F.$
Moreover, as it is known from numerical analysis, if $F'=f$ then the following equality is satisfied
\begin{equation}
\label{terer1}
f(x)=\frac{F(x+h)-F(x-h)}{2h}-\frac{h^2}{6}f''(\xi)
\end{equation}
for some $\xi\in(x-h,x+h).$ This means that \eqref{terer1} provides an alternate proof of \eqref{terHH1}
(for twice differentiable $f$).

This new formulation of the Hermite-Hadamard inequality was inspiration in \cite{terOlbrysSzostok2014} to replace
the middle term of Hermite-Hadamard inequality by more complicated expressions than those used in \eqref{terHH}. In \cite{terOlbrysSzostok2014}, the authors study inequalities of the form
$$f\left(\frac{x+y}{2}\right)\leq\frac{a_1F(x)+a_2F(\alpha x+(1-\alpha)y)+a_3F(\beta x+(1-\beta)y)+a_4F(y)}{y-x}$$
and
$$\frac{a_1F(x)+a_2F(\alpha x+(1-\alpha)y)+a_3F(\beta x+(1-\beta)y)+a_4F(y)}{y-x}\leq\frac{f(x)+f(y)}{2},$$
where $f:[x,y]\to\mathbb R$ is a convex function, $F'=f,$ $\alpha,\beta\in(0,1)$ and $a_1+a_2+a_3+a_4=0.$  
\begin{proposition}[\cite{terOlbrysSzostok2014}]
Let $n\in\mathbb{N},$ $\alpha_i\in(0,1)$, $a_i\in\mathbb R$, $i=1,\ldots,n$ be such that $\alpha_1>\alpha_2>\cdots>\alpha_n$ and $a_1+a_2+\cdots+a_n=0$, and let $F$ be a differentiable function with $F'=f.$ Then 
$$\frac{\sum_{i=1}^{n}a_iF(\alpha_i x+(1-\alpha_i)y)}{y-x}=\int fd\mu,$$
with
$$\mu(A)=-\frac{1}{y-x}\sum_{i=1}^{n-1}(a_1+\cdots+a_i)l_1(A\cap[\alpha_{i} x+(1-\alpha_{i})y,\alpha_{i+1} x+(1-\alpha_{i+1})y]),$$
where $l_1$ stands for the one-dimensional Lebesgue measure.
\label{terl1}
\end{proposition}

\begin{remark}[\cite{terOlbrysSzostok2014}]
Taking $F_1(t):=\mu((-\infty,t])$ with $\mu$ from Proposition \ref{terl1} we can see that 
\label{terR1} 
\begin{equation}
\label{terFint}
\frac{\sum_{i=1}^{n}a_iF(\alpha_i x+(1-\alpha_i)y)}{y-x}=\int fdF_1.
\end{equation}
\end{remark}

Next proposition will show that, in order to get some inequalities of the
Hermite-Hadamard type, we have to use sums containing more than three summands.
\begin{proposition}[\cite{terOlbrysSzostok2014}] \label{terl1bis}
There are no numbers $\alpha_i,a_i\in\mathbb R,i=1,2,3$, satisfying $1=\alpha_1>\alpha_2>\alpha_3=0$ such that any of the inequalities
$$f\left(\frac{x+y}{2}\right)\leq\frac{\sum_{i=1}^{3}a_iF(\alpha_i x+(1-\alpha_i)y)}{y-x}$$
or 
$$\frac{\sum_{i=1}^{3}a_iF(\alpha_i x+(1-\alpha_i)y)}{y-x}\leq\frac{f(x)+f(y)}{2}$$
is fulfilled by every continuous and convex function $f$ and its antiderivative $F.$
\end{proposition}
To prove Proposition \ref{terl1bis}, we note that by Proposition \ref{terl1}, we can see that 
$$\frac{\sum_{i=1}^{3}a_iF(\alpha_i x+(1-\alpha_i)y)}{y-x}=\int_x^yfd\mu,$$
with $$\mu(A)=-\frac{1}{y-x}\bigl(a_1l_1(A\cap[x,\alpha_2x+(1-\alpha_2)y])+$$$$
(a_2+a_1)l_1(A\cap[\alpha_2x+(1-\alpha_2)y,y])\bigr),$$ 
and
$$\frac{\sum_{i=1}^{3}a_iF(\alpha_i x+(1-\alpha_i)y)}{y-x}=\int_x^yf(t)dF_1(t),$$
where
\begin{equation}
\label{terF1}
F_1(t)=\mu\{(-\infty,t]\}.
\end{equation}
 Now, if 
$$F_2(t)=\frac{1}{y-x}l_1\{(-\infty,t]\cap[x,y]\},$$ then
$F_1$ lies strictly above or below $F_2$ (on $[x,y]$). This means that 
\begin{equation}
\label{2}
\int_{x}^yF_2(t)dt\neq\int_{x}^yF_1(t)dt.
\end{equation}
But, on the other hand, if 
\begin{equation}
\label{terF3OSZ}
F_3(t):=\left\{\begin{array}{ll}
0,\quad&t<x,\\
\frac{1}{2},\ &t\in[x,y),\\
1,\ &t\geq y,
\end{array}\right.
\end{equation}
and 
\begin{equation}
\label{terF4}
F_4(t):=\left\{\begin{array}{ll}
0,\quad & t<\frac{x+y}{2},\\
1,\quad & t\geq\frac{x+y}{2},
\end{array}\right.
\end{equation}
then 
$$\int_{x}^yF_2(t)dt=\int_{x}^yF_3(t)dt=\int_{x}^yF_4(t)dt=\frac{y-x}{2}.$$
This, together with \eqref{2}, shows that neither 
$$\int_{x}^yfdF_2\leq\int_{x}^yfdF_3$$ 
nor
$$\int_{x}^yfdF_2\geq\int_{x}^yfdF_4$$ 
is satisfied. To complete the proof it suffices to observe that 
$$\int_{x}^yfdF_3=\frac{f(x)+f(y)}{2},\;$$
$$\int_{x}^yfdF_4=f\left(\frac{x+y}{2}\right).$$
\begin{remark}[\cite{terOlbrysSzostok2014}]
 Observe that the assumptions of Proposition \ref{terl1bis}, $\alpha_1=1$ and $\alpha_3=0$, are essential.
For example, it follows from the Ohlin lemma that the inequality 
$$f\left(\frac{x+y}{2}\right)\leq\frac{-3F(\frac34 x+\frac 14 y)+\frac{25}{11}F(\frac{11}{20}x+\frac{9}{20}y)+\frac{8}{11}F(y)}{y-x}\leq\frac{1}{y-x}\int f(t)dt$$ 
is satisfied by all continuous and convex functions $f$ (where $F'=f$). Clearly, there are many 
more examples of inequalities of this type.
\end{remark}

\begin{lemma} [\cite{terOlbrysSzostok2014}]
If any of the inequalities
\begin{equation}
\label{terlhh}
f\left(\frac{x+y}{2}\right)\leq\frac{\sum_{i=1}^{4}a_iF(\alpha_i x+(1-\alpha_i)y)}{y-x}
\end{equation}
or 
\begin{equation}
\label{terrhh}
\frac{\sum_{i=1}^{4}a_iF(\alpha_i x+(1-\alpha_i)y)}{y-x}\leq\frac{f(x)+f(y)}{2}
\end{equation}
is satisfied for all continuous and convex functions $f:[x,y]\to\mathbb R$ (where $F'=f$), then 
\begin{equation}
\label{terM1}
a_1(\alpha_2-\alpha_1)+(a_2+a_1)(\alpha_3-\alpha_2)+(a_3+a_2+a_1)(\alpha_4-\alpha_3)=1
\end{equation}
and 
\begin{equation}
\label{terE1}
a_1(\alpha_2^2-\alpha_1^2)+(a_2+a_1)(\alpha_3^2-\alpha_2^2)+(a_3+a_2+a_1)(\alpha_4^2-\alpha_3^2)=1.
\end{equation}
\end{lemma}
\par\bigskip

To prove this lemma, we take $x=0$, $y=1$. Then, using Proposition \ref{terl1}, we can see that 
$$\sum_{i=1}^{4}a_iF(1-\alpha_i)=\int_0^1fd\mu=-a_1\int_{1-\alpha_1}^{1-\alpha_2}f(x)dx+$$
$$-(a_1+a_2)\int_{1-\alpha_3}^{1-\alpha_2}f(x)dx-(a_1+a_2+a_3)\int_{1-\alpha_4}^{1-\alpha_3}f(x)dx.$$
Now, we consider the functions $F_1,F_3$ and $F_4$ given by the formulas \eqref{terF1}, \eqref{terF3OSZ} and \eqref{terF4}, respectively. Then, the inequalities
\eqref{terlhh} and \eqref{terrhh} may be written in the form
$$\int fdF_4\leq\int fdF_1$$
and
$$\int fdF_1\leq\int fdF_3.$$
This means that, if for example, the inequality \eqref{terlhh} is satisfied, then we have
$F_1(1)=F_4(1)=1$, which yields \eqref{terM1}. Further, 
$$\int_{0}^1F_1(t)dt=\int_{0}^1F_4(t)dt=\frac12,$$
which gives us \eqref{terE1}.

\begin{proposition}[\cite{terOlbrysSzostok2014}]
\label{terex1}
Let $\alpha_i\in(0,1)$, $a_i\in\mathbb R$, $i=1,\dots,4,$ be
such that $1=\alpha_{1}>\alpha_2>\alpha_3>\alpha_4=0$,  $a_1+a_2+a_3+a_4=0$ and the equalities \eqref{terM1} and \eqref{terE1} are satisfied.
If $F_1$ is such that 
$$\frac{\sum_{i=1}^{4}a_iF(\alpha_{i}x+(1-\alpha_i)y)}{y-x}=\int_x^y fdF_1$$
and $F_2$ is the distribution function of a measure which is uniformly distributed in the interval $[x,y]$, then 
$(F_1,F_2)$ crosses exactly once.
\end{proposition}

Indeed, from \eqref{terM1} we can see that $F_1(x)=F_2(x)=0$ and $F_1(y)=F_2(y)=1.$ Note that,
in view of Proposition \ref{terl1}, the graph of the restriction of $F_1$ to the interval $[x,y]$ consists of three segments. Therefore, $F_1$ and $F_2$ cannot have more than one crossing point. On the other hand, if graphs $F_1$ and $F_2$ do not cross then 
$$\int_x^y tdF_1(t)\neq\int_x^y tdF_1(t)$$
i.e. \eqref{terE1} is not satisfied.
\begin{theorem}Let $\alpha_i\in(0,1)$, $a_i\in\mathbb R$, $i=1,\dots,4,$ be
such that $1=\alpha_{1}>\alpha_2>\alpha_3>\alpha_4=0$,  $a_1+a_2+a_3+a_4=0$ and the equalities \eqref{terM1} and \eqref{terE1} are satisfied.
Let $F,f:[x,y]\to\mathbb R$ be functions such that $f$ is continuous and convex and $F'=f.$ Then

(i) if $a_1>-1$, then
$$\frac{\sum_{i=1}^{4}a_iF(\alpha_ix+(1-\alpha_i)y)}{y-x}\leq\frac{1}{y-x}\int_x^yf(t)dt\leq\frac{f(x)+f(y)}{2},$$

(ii) if $a_1<-1$, then
$$f\left(\frac{x+y}{2}\right)\leq\frac{1}{y-x}\int_x^yf(t)dt\leq\frac{\sum_{i=1}^{4}a_iF(\alpha_ix+(1-\alpha_i)y)}{y-x},$$

(iii) if $a_1\in(-1,0]$, then
$$f\left(\frac{x+y}{2}\right)\leq\frac{\sum_{i=1}^{4}a_iF(\alpha_ix+(1-\alpha_i)y)}{y-x}\leq\frac{1}{y-x}\int_x^yf(t)dt,$$

(iv) if $a_1<-1$ and $a_2+a_1\leq 0$, then
$$\frac{1}{y-x}\int_x^yf(t)dt\leq\frac{\sum_{i=1}^{4}a_iF(\alpha_ix+(1-\alpha_i)y)}{y-x}\leq\frac{f(x)+f(y)}{2}.$$
\end{theorem}

We shall prove the first assertion. Other proofs are similar and will be omitted.  
It is easy to see that if inequalities which we consider are satisfied by every continuous and convex function
defined on the interval $[0,1]$, then they are true for every continuous and convex function on a given interval $[x,y].$ Therefore we assume that $x=0$ and $y=1.$ 
Let $F_1$ be such that \eqref{terFint} is satisfied and let $F_2$ be the distribution function of a measure, which is uniformly distributed in the interval $[0,1].$  From Proposition \ref{terl1} and Remark 
\ref{terR1}, we can see that the graph of $F_1$ consists of three segments and, since $a_1>-1,$
the slope of the first segment is smaller than $1,$ i.e. $F_1$ lies below $F_2$ on some 
right-hand neighborhood of $x.$  In view of the Proposition \ref{terex1}, this means that the assumptions of  the Ohlin lemma are satisfied and we get our result from this lemma. 
\par\bigskip
Now we shall present examples of inequalities, which may be obtained from this theorem.
\begin{example}[\cite{terOlbrysSzostok2014}]
\label{terex0}
Using (i), we can see that the inequality 
$${\frac13F(x)-\frac83F\left(\frac{3x+y}{4}\right)+\frac83F\left(\frac{x+3y}{4}
\right)-\frac13F(y)}\leq\frac{\int_x^yf(t)dt}{y-x}$$
is satisfied for every continuous and convex $f$ and its antiderivative $F.$
\end{example}

\begin{example}[\cite{terOlbrysSzostok2014}]
\label{terex2}
Using (ii), we can see that the inequality 
$${-2F(x)+3F\left(\frac{2x+y}{3}\right)-3F\left(\frac{x+2y}{3}
\right)+2F(y)}\geq\frac{\int_x^yf(t)dt}{y-x}$$
is satisfied by every continuous and convex function $f$ and its antiderivative $F.$
\end{example}

\begin{example}[\cite{terOlbrysSzostok2014}]
Using (iii), we can see that the inequality 
$$\frac{\int_x^yf(t)dt}{y-x}\geq\frac{-\frac12F(x)-\frac32F\left(\frac{2x+y}{3}\right)+\frac32F\left(\frac{x+2y}{3}
\right)+\frac12F(y)}{y-x}\geq f\left(\frac{x+y}{2}\right)$$
is satisfied by every continuous and convex function $f$ and its antiderivative $F.$
\end{example}

\begin{example}[\cite{terOlbrysSzostok2014}]
Using (iv), we can see that the inequality 
$$\frac{\int_x^yf(t)dt}{y-x}\leq\frac{-\frac32F(x)+2F\left(\frac{3x+y}{4}\right)-2F\left(\frac{x+3y}{4}
\right)+\frac32F(y)}{y-x}\leq\frac{f(x)+f(y)}{2}$$
is satisfied by every continuous and convex function $f$ and its antiderivative $F.$
\end{example}

In all cases considered in the above theorem, we used only the Ohlin lemma. Using Lemma \ref{terlemma:1.5}, it is possible to 
obtain more subtle inequalities. However (for the sake of simplicity), in the next result, we shall restrict our considerations to expressions of the simplified form. Note, that 
the inequality between  $f\left(\frac{x+y}{2}\right)$ and 
expressions which we consider is a bit unexpected.
\begin{theorem}[\cite{terOlbrysSzostok2014}]
\label{tert2}
Let $\alpha\in\left(0,\frac12\right)$, $a,b\in\mathbb R$.
   

(i) If $a>0$, then the inequality 
$$f\left(\frac{x+y}{2}\right)\geq\frac{aF(x)+bF(\alpha x+(1-\alpha)y)-bF((1-\alpha) x+\alpha y)-aF(y)}{y-x}$$
is satisfied by every continuous and convex $f$ and its antiderivative $F$ if, and only if,
\begin{equation}
\label{terA1A2}
(1-\alpha)^2\frac{ab}{a+b}>\frac12-(1-\alpha)\frac{b}{a+b},
\end{equation}

(ii) if $a<-1$ and $a_1+a_2>0$, then the inequality 
\begin{equation*}
\frac{aF(x)+bF(\alpha x+(1-\alpha)y)-bF((1-\alpha) x+\alpha y)-aF(y)}{y-x}\leq \frac{f(x)+f(y)}{2}
\end{equation*}
is satisfied by every continuous and convex $f$ and its antiderivative $F$ if, and only if,
$$-\frac{1}{4a}>\left(-a(1-\alpha)-\frac12\right)\left(\frac12+\frac{1}{2a}\right).$$
\end{theorem}
We shall prove the assertion (i) of Theorem \ref{tert2}. The proof of (ii) is similar and will be omitted. 
Similarly as before, we may assume without loss of generality, that $x=0,y=1$. Let $F_1$ be such that 
$$aF(0)+bF(1-\alpha)-bF(\alpha)+aF(1)=\int_0^1fdF_1$$
and let $F_4$ be given by \eqref{terF4}. Then it is easy to see that $(F_1,F_4)$ crosses three times:
at $\frac{(1-\alpha)b}{a+b},$ $\frac12$ and at $\frac{a+\alpha b}{a+b},$

We are going to use Lemma \ref{terlemma:1.5}. Since, from \eqref{terE1}, we have that
$$A_0+A_1+A_2+A_3=0,$$ 
it suffices to check that 
$A_0\geq A_1$ if, and only if, the inequality \eqref{terA1A2} is satisfied.
Since, $F_4(x)=0,$ for  $x\in\left(0,\frac12\right),$ we get 
$$A_0=-\int_{0}^{\frac{(1-\alpha)b}{a+b}}F_1(t)dt$$
and 
$$A_1=\int_{\frac{(1-\alpha)b}{a+b}}^\frac12F_1(t)dt,$$
which yields our assertion.
\begin{example}[\cite{terOlbrysSzostok2014}]
Neither inequality
\begin{equation}
\label{ter*}
f\left(\frac{x+y}{2}\right)\leq\frac{\frac13F(x)-\frac83F\left(\frac{3x+y}{4}\right)+\frac83F\left(\frac{x+3y}{4}
\right)-\frac13F(y)}{y-x}
\end{equation}
 nor 
 \begin{equation}
 \label{ter**}
f\left(\frac{x+y}{2}\right)\geq\frac{\frac13F(x)-\frac83F\left(\frac{3x+y}{4}\right)+\frac83F\left(\frac{x+3y}{4}
\right)-\frac13F(y)}{y-x}
\end{equation}
is satisfied for all continuous and convex $f:[x,y]\to\mathbb R.$ Indeed, if $F_1$ is such that  
$$\int_x^yf(t)dF_1(t)=\frac{\frac13F(x)-\frac83F\left(\frac{3x+y}{4}\right)+\frac83F\left(\frac{x+3y}{4}
\right)-\frac13F(y)}{y-x},$$
then 
$$\int_x^{\frac{3x+y}{4}}F_1(t)dt<\int_x^{\frac{3x+y}{4}}F_4(t)dt,$$
thus inequality \eqref{ter*} cannot be satisfied. On the other hand, the coefficients and nodes 
of the expression considered do not satisfy \eqref{terA1A2}. Therefore \eqref{ter**} is also not satisfied 
for all continuous and convex $f:[x,y]\to\mathbb R.$
\end{example}

\begin{example}[\cite{terOlbrysSzostok2014}]
Using assertion (i) of Theorem \ref{tert2}, we can see that the inequality 
$$\frac{2F(x)-3F\left(\frac{3x+y}{4}\right)+3F\left(\frac{x+3y}{4}
\right)-2F(y)}{y-x}\leq f\left(\frac{x+y}{2}\right)$$
is satisfied for every continuous and convex $f$ and its antiderivative $F.$
\end{example}

\begin{example}[\cite{terOlbrysSzostok2014}]
Using  assertion (ii) of Theorem \ref{tert2}, we can see that the inequality 
$$\frac{-2F(x)+3F\left(\frac{2x+y}{3}\right)-3F\left(\frac{x+2y}{3}
\right)+2F(y)}{y-x}\leq \frac{f(x)+f(y)}{2}$$
is satisfied for every continuous and convex $f$ and its antiderivative $F.$
\end{example}
\section{Inequalities of the Hermite-Hadamard type involving numerical differentiation formulas of  order two}
\label{terdifferentiation two}
In the paper \cite{terTSzostok2016b}, expressions connected with numerical differentiation formulas of order 
$2$, are studied. The author used the Ohlin lemma and the Levin-Ste\v{c}kin
theorem to study inequalities connected with these expressions.
In particular, the author present a new proof of the inequality 
\begin{equation}
\label{terDr}
 f\left(\frac{x+y}{2}\right)\leq\frac{1}{(y-x)^2}\int_x^y\hspace{-2mm}\int_x^yf\left(\frac{s+t}{2}\right)ds\:dt
\leq\frac{1}{y-x}\int_x^yf(t)dt,
\end{equation}
satisfied by every convex function $f:\mathbb R\to\mathbb R$ 
and he obtain extensions of 
  \eqref{terDr}. 
 In the previous section, inequalities involving expressions of the form 
$$\frac{\sum_{i=1}^na_iF(\alpha_ix+\beta_iy)}{y-x},$$
where $\sum_{i=1}^na_i=0,$ $\alpha_i+\beta_i=1$ and $F'=f$
were considered. 
 In this section, we study inequalities 
for expressions of the form 
\begin{equation*}
\frac{\sum_{i=1}^na_iF(\alpha_ix+\beta_iy)}{(y-x)^2},
\end{equation*}
which we use to approximate the second order derivative of $F$ and, surprisingly,
we discover a connection between our approach and the inequality \eqref{terDr} (see \cite{terTSzostok2016b}).

 First, we make the following simple observation.
\begin{remark}[\cite{terTSzostok2016b}]
\label{01}
Let $f,F,\Phi:[x,y]\to\mathbb R$ be such that $\Phi'=F,F'=f$. Let $n_i,m_i\in\mathbb N\cup\{0\}$, $i=1,2,3$; $a_{i,j}\in\mathbb R$, $\alpha_{i,j}$, $\beta_{i,j}\in[0,1]$, $i=1,2,3$; $j=
1,\dots,n_i,$
$b_{i,j}\in\mathbb R$, $\gamma_{i,j}$, $\delta_{i,j}\in[0,1]$, $i=1,2,3$; $j=1,\dots,m_i.$
If the inequality
\begin{multline}
\label{ter01e}
\sum_{i=1}^{n_1}a_{1,i}f(\alpha_{1,i}x+\beta_{1,i}y)+\frac{\sum_{i=1}^{n_2}a_{2,i}F(\alpha_{2,i}x+\beta_{2,i}y)}{y-x}\\+\frac{\sum_{i=1}^{n_3}a_{3,i}\Phi(\alpha_{3,i}x+\beta_{3,i}y)}{(y-x)^2}
\leq\sum_{i=1}^{m_1}b_{1,i}f(\gamma_{1,i}x+\delta_{1,i}y)\\+\frac{\sum_{i=1}^{m_2}b_{2,i}F(\gamma_{2,i}x+\delta_{2,i}y)}{y-x}+\frac{\sum_{i=1}^{m_3}b_{3,i}\Phi(\gamma_{3,i}x+\delta_{3,i}y)}{(y-x)^2}
\end{multline}
is satisfied for $x=0,y=1$ and for all continuous and convex functions $f:[0,1]\to\mathbb R$, then it is satisfied  
for all $x,y\in\mathbb R$, $x<y$ and for each continuous and convex function $f:[x,y]\to\mathbb R.$
To see this it is enough to observe that  expressions from \eqref{ter01e} remain unchanged if we replace $f:[x,y]\to\mathbb R$ 
by $\varphi:[0,1]\to\mathbb R$ given by $\varphi(t):=f\left(x+\frac{t}{y-x}\right).$
\end{remark}

The simplest expression used to approximate the second order derivative of $f$ is of the form 
$$f''\left(\frac{x+y}{2}\right)\approx\frac{f(x)-2f\left(\frac{x+y}{2}\right)+f(y)}{\left(\frac{y-x}{2}\right)^2}$$
\begin{remark}[\cite{terTSzostok2016b}]
\label{ternuman}
From numerical analysis it is known that 
$$f''\left(\frac{x+y}{2}\right)=\frac{f(x)-2f\left(\frac{x+y}{2}\right)+f(y)}{\left(\frac{y-x}{2}\right)^2}-
\frac{\left(\frac{y-x}{2}\right)^2}{12}f^{(4)}(\xi).$$
This means that for a convex function $g$ and for $G$ such that $G''=g$ we have
$$g\left(\frac{x+y}{2}\right)\leq\frac{G(x)-2G\left(\frac{x+y}{2}\right)+G(y)}{\left(\frac{y-x}{2}\right)^2}.$$
\end{remark}
In the paper \cite{terTSzostok2016b}, some inequalities for convex functions which do not follow from 
formulas used in numerical differentiation, are obtained .

Let now $f:[x,y]\to\mathbb R$ be any function and let $F,\Phi:[x,y]\to\mathbb R$ be such that $F'=f$ and $\Phi''=f.$
We need to write the expression 
\begin{equation}
\label{terP}
\frac{\Phi(x)-2\Phi\left(\frac{x+y}{2}\right)+\Phi(y)}{\left(\frac{y-x}{2}\right)^2}
\end{equation}
in the form $$\int_x^yfdF_1$$ for some $F_1.$
In the next proposition we show that it is possible -- here for the sake of simplicity we shall work on the interval $[0,1].$
\begin{proposition}[\cite{terTSzostok2016b}]
Let $f:[0,1]\to\mathbb R$ be any function and let $\Phi:[0,1]\to\mathbb R$ be such that $\Phi''=f.$ Then we have 
$$4\left(\Phi(0)-2\Phi\left(\frac{1}{2}\right)+\Phi(1)\right)=\int_x^yfdF_1,$$
where $F_1:[0,1]\to\mathbb R$ is given by
\begin{equation}
\label{terF3}
F_1(t):=\left\{\begin{array}{ll}
2x^2,\ &x\leq\frac12,\\
-2x^2+4x-1,\ &x>\frac12.
\end{array}\right.
\end{equation}
\end{proposition}
Now, we observe that the following equality is satisfied
$$\frac{\Phi(x)-2\Phi\left(\frac{x+y}{2}\right)+\Phi(y)}{\left(\frac{y-x}{2}\right)^2}=\frac{1}{(y-x)^2}
\int_x^y\hspace{-2mm}\int_x^y
f\left(\frac{s+t}{2}\right)ds\:dt.$$
After this observation it turns out that inequalities involving the expression \eqref{terP} were considered 
in the paper of Dragomir \cite{terDragomir}, where (among others) the following inequalities were obtained
\begin{equation}
\label{terD}
 f\left(\frac{x+y}{2}\right)\leq\frac{1}{(y-x)^2}\int_x^y\hspace{-2mm}\int_x^yf\left(\frac{s+t}{2}\right)ds\:dt
\leq\frac{1}{y-x}\int_x^yf(t)dt.
\end{equation}
As we already know (Remark \ref{ternuman}) the first one of the above inequalities may be obtained using the numerical analysis
results. 

Now, the inequalities from the Dragomir's paper easily follow from the Ohlin lemma but 
there are many possibilities of generalizations and modifications of inequalities 
\eqref{terD}. These generalizations will be discussed in this section.
\par\bigskip

First, we consider the symmetric case. We start with the following remark. 
\begin{remark}[\cite{terTSzostok2016b}]
\label{terrma}
Let $F_*(t)=at^2+bt+c$ for some $a,b,c\in\mathbb R, a\neq 0.$
It is impossible to obtain inequalities involving
$\int_x^yfdF_*$ and any of the expressions: 
$$\frac{1}{y-x}\int_x^yf(t)dt,\quad \ \ f\left(\frac{x+y}{2}\right),\quad \ \ \frac{f(x)+f(y)}{2},$$
which are satisfied for all convex functions $f:[x,y]\to\mathbb R.$ Indeed, suppose that we have 
$$\int_x^yfdF_*\leq\frac{1}{y-x}\int_x^yf(t)dt$$
for all convex $f:[x,y]\to\mathbb R.$
Without loss of generality we may assume that $F_*(x)=0,$ then from Theorem \ref{terth:1.3} we have $F_*(y)=1$.
Also from Theorem \ref{terth:1.3} we get
$$\int_x^yF_*(t)dt=\int_x^yF_0dt,$$
where $F_0(t)=\frac{t-x}{y-x}$, $t\in[x,y]$,
which is impossible, because $F_*$ is either strictly convex or concave.
 \end{remark}

\par\bigskip
 This remark means that in order to get some new inequalities of the Hermite-Hadamard type we have 
to integrate with respect to functions constructed with use of (at least) two quadratic functions.

Now we present the main result of this section.
\begin{theorem}[\cite{terTSzostok2016b}]
\label{termainth}
Let $x,y$ be some real numbers such that $x<y$ and let $a\in\mathbb R.$
Let $f,F,\Phi:[x,y]\to\mathbb R$ be any functions such that $F'=f$ and $\Phi'=F$ and let $T_af(x,y)$ be the function defined
by the following formula 
$$T_af(x,y)=\left(1-\frac{a}{2}\right)\frac{F(y)-F(x)}{y-x}+2a\frac{\Phi(x)-2\Phi\left(\frac{x+y}{2}\right)+\Phi(x)}{(y-x)^2}.$$
Then the following inequalities hold for all convex functions $f:[x,y]\to\mathbb R:$
\begin{itemize}
\item if $a\geq 0$, then
\begin{equation}
\label{terai}
T_af(x,y)\leq\frac{1}{y-x}\int_x^yf(t)dt,
\end{equation}   
\item if $a\leq 0$, then
\begin{equation}
\label{teraii}
T_af(x,y)\geq\frac{1}{y-x}\int_x^yf(t)dt,
\end{equation}   
\item if $a\leq 2$, then
\begin{equation}
\label{teraiii}
f\left(\frac{x+y}{2}\right)\leq T_af(x,y),
\end{equation}   
\item if $a\geq 6$, then
\begin{equation}
\label{teraiv}
T_af(x,y)\leq f\left(\frac{x+y}{2}\right),
\end{equation}   
\item if $a\geq -6$, then
\begin{equation}
\label{terav}
T_af(x,y)\leq\frac{f(x)+f(y)}{2},
\end{equation}
\end{itemize}   
Furthermore,
\begin{itemize}
 \item if $a\in(2,6)$, then the expressions $T_af(x,y),$ $f\left(\frac{x+y}{2}\right)$ 
are not comparable in the class of convex functions,

\item  if $a<-6$, then
expressions $T_af(x,y),$ $\frac{f(x)+f(y)}{2}$ are not comparable in the class of convex functions.
\end{itemize}
\end{theorem}
To prove Theorem \ref{termainth}, we note that, we may restrict ourselves to the case $x=0,y=1.$
Take $a\in\mathbb R,$ let $f:[0,1]:\to\mathbb R$ be any convex function and let $F,\Phi:[0,1]\to\mathbb R$ be such that $F'=f,\Phi'=F.$  
Define  $F_1:[0,1]\to\mathbb R$ by the formula
\begin{equation}
\label{terF1t}
F_1(t):=\left\{\begin{array}{ll}
at^2+\left(1-\frac{a}{2}\right)t,\quad &t<\frac12,\\
-at^2+\left(1+\frac{3a}{2}\right)t-\frac{a}{2},\quad &t\geq\frac12.
\end{array}\right.
\end{equation}
First, we prove that $T_af(0,1)=\int_0^1fdF_1.$
Now, let $F_2(t)=t$, $t\in[0,1].$
Then the functions $F_1,F_2$ have exactly one crossing point (at $\frac12$)
and 
$$\int_0^1F_1(t)dt=\frac12=\int_0^1tdt.$$
Moreover, if $a>0$, then the function $F_1$ is convex on the interval $(0,\frac12)$ and 
concave on $(\frac12,1).$ 
Therefore, it follows from the Ohlin lemma, that
for $a>0$ we have
\begin{equation*}
\int_0^1fdF_1\leq\int_0^1fdF_2,
\end{equation*}
which, in view of Remark \ref{01}, yields \eqref{terai} and for $a<0$ the opposite inequality is satisfied, which gives \eqref{teraii}.
Take
\begin{equation*}
F_3(t):=\left\{\begin{array}{ll}
0,\quad &t\leq\frac12,\\
1,\quad &t>\frac12.
\end{array}\right.
\end{equation*}
 
It is easy to calculate that for $a\leq 2$ we have $F_1(t)\geq F_3(t)$ for $t\in\left[0,\frac 12\right],$
and $F_1(t)\leq F_3(t)$ for $t\in\left[\frac 12,1\right]$, and this means that from 
the Ohlin lemma we get \eqref{teraiii}.
Let now
\begin{equation*}
F_4(t):=\left\{\begin{array}{ll}
0,\quad &t=0,\\
\frac12,\quad &t\in(0,1),\\
1,\quad &t=1.
\end{array}\right.
\end{equation*}
Similarly as before, if $a\geq-2$, then we have  $F_1(t)\geq F_4(t)$ for $t\in\left[0,\frac 12\right]$ and 
$F_1(t)\leq F_4(t)$ for $t\in\left[\frac 12,1\right].$ Therefore, from the Ohlin lemma, we get \eqref{teraiv}. 

Suppose that $a>2.$ Then there are three crossing points of the functions $F_1$ and $F_3$ $:x_0,\frac12,x_1,$
where $x_0\in(0,\frac12), x_1\in(\frac12,1)$. The function 
$$\varphi(s):=\int_0^s(F_3(t)-F_1(t))dt,\;s\in[0,1]$$
is increasing on the intervals $[0,x_0],[\frac12,x_1]$ and decreasing on $[x_0,\frac12]$ and on $[x_1,1].$ 
This means that $\varphi$ takes its absolute minimum at $\frac12.$
It is easy to calculate that $\varphi\left(\frac12\right)\geq 0$, if $a\geq 6$, which, in view of  Theorem \ref{terth:1.3},
gives us \eqref{teraiv}.

To see, that for $a\in(2,6)$, the expressions $T_af(x,y)$ and $f\left(\frac{x+y}{2}\right)$ are not comparable in the class 
of convex functions it is enough to observe that in this case $\varphi(x_0)>0$ and $\varphi\left(\frac12\right)<0.$

Analogously (using functions $F_1$ and $F_4$), we show that for $a\in(-2,-6]$ we have \eqref{terav}, and in the case $a<-6$ the expressions
$T_af(x,y)$ and $\frac{f(x)+f(y)}{2}$ are not comparable in the class of convex functions.
This theorem provides us with a full description of inequalities, which may be obtained using  Stieltjes integral
with respect to a function of the form \eqref{terF1t}. Some of the obtained inequalities are already known. For example, from 
\eqref{terai} and \eqref{teraii} we obtain the inequality 
$$\frac{1}{(y-x)^2}\int_x^y\hspace{-2mm}\int_x^yf\left(\frac{s+t}{2}\right)ds\:dt
\leq\frac{1}{y-x}\int_x^yf(t)dt,$$
whereas from \eqref{teraiii} for $a=2$ we get the inequality 
$$f\left(\frac{x+y}{2}\right)\leq\frac{1}{(y-x)^2}\int_x^y\hspace{-2mm}\int_x^yf\left(\frac{s+t}{2}\right)ds\:dt.$$

 However, inequalities obtained  for "critical" values of $a$
i.e. $-6,6.$ are here particularly interesting. 
 In the following corollary, we explicitly write these inequalities.
\begin{corollary}[\cite{terTSzostok2016b}]
\label{terwn}
For every convex function $f:[x,y]\to\mathbb R$, the following inequalities are satisfied

\begin{equation}
\label{teri}
3\frac{1}{(y-x)^2}\int_x^y\int_x^yf\left(\frac{s+t}{2}\right)dsdt
\leq \frac{2}{y-x}\int_x^yf(t)dt+f\left(\frac{x+y}{2}\right),
\end{equation}
\begin{equation}
\label{terii}
\frac{4}{y-x}\int_x^yf(t)dt\leq 3\frac{1}{(y-x)^2}\int_x^y\int_x^yf\left(\frac{s+t}{2}\right)dsdt+\frac{f(x)+f(y)}{2}.
\end{equation}
 \end{corollary}
\begin{remark}[\cite{terTSzostok2016b}]
In the paper \cite{terDG}, S.S. Dragomir and I. Gomm obtained the following inequality 
\begin{equation}
\label{terdrgo}
3\int_x^yf(t)dt\leq 2\frac{1}{(y-x)^2}\int_x^y\int_x^yf\left(\frac{s+t}{2}\right)dsdt+\frac{f(x)+f(y)}{2}.
\end{equation}
Inequality \eqref{terii} from Corollary \ref{terwn} is stronger than \eqref{terdrgo}. Moreover, as it was observed in Theorem \ref{termainth},
the inequalities \eqref{teri} and \eqref{terii} cannot be improved i.e. the inequality
$$\frac{1}{y-x}\int_x^yf(t)dt\leq \lambda\frac{1}{(y-x)^2}\int_x^y\int_x^yf\left(\frac{s+t}{2}\right)dsdt+(1-\lambda)\frac{f(x)+f(y)}{2}$$
for $\lambda>\frac34$ is not satisfied by every convex function $f:[x,y]\to\mathbb R$
and the inequality 
$$\frac{1}{(y-x)^2}\int_x^y\int_x^yf\left(\frac{s+t}{2}\right)dsdt
\leq \gamma\frac{1}{y-x}\int_x^yf(t)dt+(1-\gamma)f\left(\frac{x+y}{2}\right)
$$
with $\gamma>\frac23$ is not true for all convex functions  $f:[x,y]\to\mathbb R.$ 
\end{remark} 
In Corollary \ref{terwn} we obtained inequalities for the triples:
$$\frac{1}{(y-x)^2}\int_x^y\hspace{-1.5mm}\int_x^yf\left(\frac{s+t}{2}\right)dsdt,\quad \int_x^yf(t)dt,\quad \frac{f(x)+f(y)}{2}$$
and
$$\frac{1}{(y-x)^2}\int_x^y\hspace{-1.5mm}\int_x^yf\left(\frac{s+t}{2}\right)dsdt,\quad \int_x^yf(t)dt,\quad f\left(\frac{x+y}{2}\right).$$
In the next remark, we present an analogous result for expressions
$$\frac{1}{(y-x)^2}\int_x^y\hspace{-1.5mm}\int_x^yf\left(\frac{s+t}{2}\right)dsdt,\quad \frac{f(x)+f(y)}{2},\quad f\left(\frac{x+y}{2}\right).$$
\begin{remark}[\cite{terTSzostok2016b}]
Using the functions: $F_1$ defined by \eqref{terF3} and $F_5$ given by 
\begin{equation*}
\label{terF5t}
F_5(t):=\left\{\begin{array}{ll}
0,& \quad t=0,\\
\frac16,&\quad t\in\left(0,\frac12\right),\\
\frac56,&\quad t\in\left[\frac12,1\right),\\
1,&\quad t=1,
\end{array}\right.
\end{equation*}
we can see that 
$$\frac16f(x)+\frac23f\left(\frac{x+y}{2}\right)+\frac16f(y)\geq\frac{1}{(y-x)^2}\int_x^y\hspace{-1.5mm}\int_x^yf\left(\frac{s+t}{2}\right)dsdt$$
for all convex functions $f:[x,y]\to\mathbb R.$

Moreover, it is easy to see, that 
 the above inequality cannot be strengthened, which  means that, 
if $a$ $,b\geq 0$, $2a+b=1$ and $a<\frac16$, then
the inequality 
$$af(x)+bf\left(\frac{x+y}{2}\right)+af(y)\geq\frac{1}{(y-x)^2}\int_x^y\hspace{-1.5mm}\int_x^yf\left(\frac{s+t}{2}\right)dsdt,$$
is not satisfied by all convex functions $f$.
\end{remark}

In \cite{terTSzostok2016b}, inequalities for  $f(\alpha x+(1-\alpha)y)$ and for $\alpha f(x)+(1-\alpha) f(y),$ 
where $\alpha$ is not necessarily equal to $\frac12$ (the non-symmetric case), are also obtained.
\begin{theorem}[\cite{terTSzostok2016b}]
\label{terNS}
Let $x,y$ be some real numbers such that $x<y$ and let $\alpha\in[0,1].$ Let  $f:[x,y]\to\mathbb R,$
be a convex function, let 
 $F$ be such that $F'=f$ and let $\Phi$ satisfy $\Phi'=F.$  
If $S^2_\alpha f(x,y)$ is defined by 
$$S^2_\alpha f(x,y):=\frac{(4-6\alpha)F(y)+(2-6\alpha)F(x)}{y-x}+\frac{(6-12\alpha)(\Phi(y)-\Phi(x))}{(y-x)^2},$$
then the following conditions hold true:
\begin{itemize}
\item \begin{equation*}
\label{teral}
S^2_\alpha f(x,y)\leq\alpha f(x)+(1-\alpha) f(y),
\end{equation*}   

\item if $\alpha\in\left[\frac13,\frac23\right]$, then
\begin{equation*}
\label{teral1}
S^2_\alpha f(x,y)\geq f(\alpha x+(1-\alpha)y),
\end{equation*}    
\item if $\alpha\in[0,1]\setminus\left[\frac13,\frac23\right]$, then the expressions $S^2_\alpha f(x,y)$ and $f(\alpha x+(1-\alpha)y)$ are incomparable in 
the class of convex functions,

\item if $\alpha\in\left(0,\frac13\right]\cup\left[\frac23,1\right),$ then
\begin{equation*}
\label{terS1S2}
S^2_\alpha f(x,y)\leq S^1_\alpha f(x,y),
\end{equation*}    
\item  if $\alpha\in\left(\frac13,\frac12\right)\cup\left(\frac12,\frac23\right)$, then $S^1_\alpha f(x,y)$ and $S^2_\alpha f(x,y)$ are incomparable in 
the class of convex functions.
\end{itemize}
\end{theorem}

\section{The  Hermite-Hadamard type inequalities for $n$-th order convex functions}
\label{tern order}

Now we are going to study Hermite-Hadamard type inequalities for higher-order convex functions. Many results on higher order generalizations of the Hermite-Hadamard type inequality one can found, among others, in \cite{terBes08, terBesPal02, terBesPal03, terBesPal04, terBesPal10, terDraPe2000, terBesPal08, terRajba12f, terRajba12g}. In recent papers \cite{terRajba12f, terRajba12g}, the theorem of M. Denuit, C.Lef\`{e}vre and  M. Shaked \cite{terDenLefSha98} on sufficient conditions for $s$-convex ordering was used, to prove Hermite-Hadamard type inequalities for higher-order convex functions.

Let us review some notations. The convexity of $n$-th order (or $n$-convexity) was defined in terms of divided differences by Popoviciu \cite{terPopoviciu1934}, however, we will not state it here. Instead we list some properties of $n$-th order convexity which are equivalent to Popoviciu's definition (see \cite{terKuczma1985}).

\begin{proposition}\label{terprop:16}
A function $f \colon (a,b)\to\mathbb{R}$ is $n$-convex on $(a,b)$ $(n \geq 1)$ if, and only if, its derivative $f^{(n-1)}$ exists and is convex on $(a,b)$ (with the convention $f^{(0)}(x) = f(x)$).
\end{proposition}

\begin{proposition}\label{terprop:17}
Assume that $f \colon [a,b] \to\mathbb{R}$ is $(n+1)$-times differentiable on $(a,b)$ and continuous on $[a,b]$ ($n \geq 1$). Then $f$ is $n$-convex if, and only if, $f^{(n+1)}(x)\geq 0$, $x \in (a,b)$.
\end{proposition}

For real valued random variables $X,Y$ and any integer $s \geq 2$ we say that $X$ is dominated by $Y$ in $s$-\textit{convex ordering} sense if $\mathbb{E} f(X) \leq \mathbb{E} f(Y)$ for all $(s-1)$-convex functions $f \colon \mathbb{R} \to \mathbb{R}$, for which the expectations exist (\cite{terDenLefSha98}). In that case we write $X \leq_{s-cx} Y$, or $\mu_X \leq_{s-cx} \mu_Y$, or $F_X \leq_{s-cx} F_Y$. Then the order $\leq_{2-cx}$ is just the usual convex order $\leq_{cx}$.

A very useful criterion for the verification of the $s$-convex order is given by Denuit, Lef\`{e}vre and Shaked in \cite{terDenLefSha98}.
\begin{proposition}[\cite{terDenLefSha98}]\label{terprop:18}
Let $X$ and $Y$ be two random variables such that $\mathbb{E} (X^j-Y^j)=0$, $j=1,2,\ldots,s-1$ ($s \geq 2$). If $S^-(F_X-F_Y)=s-1$ and the last sign of $F_X-F_Y$ is positive, then $X \leq_{s-cx} Y$.
\end{proposition}

We now apply Proposition \ref{terprop:18} to obtain the following results.

\begin{theorem}[\cite{terRajba12f}]\label{terth:19}
Let $n\geq 1$, $a_1 \leq a < b \leq b_1$. 

Let $a(n) = \left[\frac{n}{2}\right]+1$, $b(n) = \left[\frac{n+1}{2}\right]+1$.

Let $\alpha_1,\ldots,\alpha_{a(n)}$, $x_1,\ldots,x_{a(n)}$, $\beta_1,\ldots,\beta_{b(n)}$, $y_1,\ldots,y_{b(n)}$ be real numbers such that
\begin{itemize}
	\item if $n$ is even then
	\begin{eqnarray*}
	& 0 < \beta_1 < \alpha_1 < \beta_1+\beta_2 < \alpha_1 + \alpha_2 < \ldots < \alpha_1+\ldots+\alpha_{a(n)} = \beta_1+\ldots+\beta_{b(n)}=1, & \nonumber \\
	& a \leq y_1 < x_1 < y_2 < x_2 < \ldots < x_{a(n)} < y_{b(n)} \leq b,
	\end{eqnarray*}
	\item if $n$ is odd then
	\begin{eqnarray*}
		& 0 < \beta_1 < \alpha_1 < \beta_1 + \beta_2 < \alpha_1+\alpha_2 < \ldots < \beta_1+\ldots+\beta_{b(n)} < \alpha_1+\ldots+\alpha_{a(n)}=1 & \nonumber \\
		& a \leq y_1 < x_1 < y_2 < x_2 < \ldots < y_{b(n)} < x_{a(n)} \leq b;
	\end{eqnarray*}
\end{itemize}
and
$$
\sum_{k=1}^{a(n)} x_i^k \alpha_i = \sum_{j=1}^{b(n)} y_j^k\beta_j
$$
for any $k=1,2,\ldots,n$.

Let $f \colon [a_1,b_1] \to \mathbb{R}$ be an $n$-convex function. Then we have the following inequalities:
\begin{itemize}
	\item if $n$ is even then

\begin{equation*}
\sum_{i=1}^{a(n)} \alpha_i f(x_i) \leq \sum_{j=1}^{b(n)} \beta_j f(y_j),
\label{eq:15}
\end{equation*}
\item if $n$ is odd then
\begin{equation*}
\sum_{j=1}^{b(n)} \beta_j f(y_j) \leq \sum_{i=1}^{a(n)} \alpha_i f(x_i).
\label{eq:16}
\end{equation*}
\end{itemize}
\end{theorem}

\begin{theorem}[\cite{terRajba12f}]\label{terth:19a}
Let $n\geq 1$, $a_1 \leq a < b \leq b_1$. Let $a(n), b(n) \in \mathbb{N}$. Let $\alpha_1,\ldots,\alpha_{a(n)}$, $\beta_1,\ldots,\beta_{b(n)}$ be positive real numbers such that  $\alpha_1+\ldots+\alpha_{a(n)} = \beta_1+\ldots+\beta_{b(n)}=1$. Let $x_1,\ldots,x_{a(n)}$,  $y_1,\ldots,y_{b(n)}$ be real numbers such that
\begin{itemize}
\item
$a \leq x_1 \leq x_2 \leq \ldots\leq x_{a(n)} \leq b$ and $a \leq y_1\leq y_2 \leq \ldots \leq y_{b(n)} \leq b$,
\item
$
\sum_{k=1}^{a(n)} x_i^k \alpha_i = \sum_{j=1}^{b(n)} y_j^k\beta_j,
$
for any $k=1,2,\ldots,n$.
\end{itemize}
Let $\alpha_0=\beta_0 = 0$, $x_0=y_0= - \infty$. Let $F_1, F_2 \colon \mathbb{R} \to \mathbb{R}$ be two functions given by the following formulas: $F_1 (x)= \alpha_0+ \alpha_1+\ldots+\alpha_k$ if $ x_k<x\leq x_{k+1}$ $(k=0,1, \ldots , a(n)-1)$ and $F_1(x)=1$ if $x>x_{a(n)}$; $F_2 (x)= \beta_0+ \beta_1+\ldots+\beta_k$ if $ y_k<x\leq y_{k+1}$ $(k=0,1, \ldots , b(n)-1)$ and $F_2(x)=1$ if $x>y_{b(n)}$. If the functions $F_1, F_2$ have $n $ crossing points and the last sign of $F_1-F_2$ is a+, then  for any $n$-convex function $f \colon [a_1,b_1] \to \mathbb{R}$ we have the following inequality
\begin{equation*}
\sum_{i=1}^{a(n)} \alpha_i f(x_i) \leq \sum_{j=1}^{b(n)} \beta_j f(y_j).
\end{equation*}
\end{theorem}

\begin{theorem}[\cite{terRajba12f}]\label{terth:20}
Let $n\geq 1$, $a_1 \leq a < b \leq b_1$. Let $a(n) = \left[\frac{n}{2}\right]+1$, $b(n) = \left[\frac{n+1}{2}\right]+1$. Let $x_1,\ldots,x_{a(n)},y_1,\ldots,y_{b(n)}$ be real numbers, and $\alpha_1,\ldots,\alpha_{a(n)}$, $\beta_1,\ldots,\beta_{b(n)}$ be  positive numbers, such that $\alpha_1+\ldots+\alpha_{a(n)}=1$, $\beta_1+\ldots+\beta_{b(n)}=1$,
$$
\frac{1}{b-a}\int_a^b x^k dx = \sum_{j=1}^{b(n)} y_j^k \beta_j = \sum_{i=1}^{a(n)} x_i^k \alpha_i \quad (k=1,2,\ldots,n),
$$
$a \leq x_1 < x_2 < \ldots < x_{a(n)} \leq b$, $a \leq y_1 < y_2 < \ldots < y_{b(n)} < b$,
\begin{eqnarray*}
& \frac{x_1-a}{b-a} < \alpha_1 < \frac{x_2-a}{b-a}, & \\
& \frac{x_2-a}{b-a} < \alpha_1+\alpha_2 < \frac{x_3-a}{b-a}, & \\
& \ldots & \\
& \frac{x_{a(n)-1}-a}{b-a} < \alpha_1+\ldots+\alpha_{a(n)-1} < \frac{x_{a(n)}-a}{b-a}, &
\end{eqnarray*}
\begin{eqnarray*}
& \frac{y_1-a}{b-a} < \beta_1 < \frac{y_2-a}{b-a}, & \\
& \frac{y_2-a}{b-a} < \beta_1+\beta_2 < \frac{y_2-a}{b-a}, & \\
& \ldots & \\
& \frac{y_{b(n)-1}-a}{b-a} < \beta_1+\ldots+\beta_{b(n)-1} < \frac{y_{b(n)}-a}{b-a}; &
\end{eqnarray*}
if $n$ is even then $y_1=a$, $y_{b(n)}=b$, $x_1>a$, $x_{a(n)}<b$; \\
if $n$ is odd then $y_1=a$, $y_{b(n)}<b$, $x_1>a$, $x_{a(n)}=b$.

Let $f \colon [a_1,b_1] \to\mathbb{R}$ be an $n$-convex function. Then we have the following inequalities:\\
\begin{itemize}
\item if $n$ is even then
\begin{equation*}
\sum_{i=1}^{a(n)} \alpha_i f(x_i) \leq \frac{1}{b-a} \int_a^b f(x)dx \leq \sum_{j=1}^{b(n)} \beta_j f(y_j),
\end{equation*}
\item if $n$ is odd then
\begin{equation*}
\sum_{j=1}^{b(n)} \beta_j f(y_j) \leq \frac{1}{b-a} \int_a^b f(x) dx \leq \sum_{i=1}^{a(n)} \alpha_i f(x_i).
\end{equation*}
\end{itemize}
\end{theorem}
Note, that Proposition \ref{terprop:18} can be rewritten in the following form.
\begin{proposition}[\cite{terDenLefSha98}]\label{terprop:19}
Let $X$ and $Y$ be two random variables such that 
\begin{equation*}
\mathbb{E} (X^j-Y^j)=0,\quad j=1,2,\ldots,s \:(s \geq 1).
\end{equation*}
If the distribution functions $F_X$ and $F_Y$ cross exactly $s$-times at points $x_1<x_2< \ldots <x_s$ and 
\begin{equation*}
(-1)^{s+1}\left(F_Y(x)-F_X(x)\right)\geq 0\quad for\: all\: x\leq x_1,
\end{equation*}
then 
\begin{equation}\label{eq:ord5}
\mathbb{E} f(X) \leq \mathbb{E} f(Y)
\end{equation}
for all $s$-convex functions $f \colon \mathbb{R} \to \mathbb{R}$.
\end{proposition}

Proposition \ref{terprop:18} is a counterpart of the Ohlin lemma concerning convex ordering.  This proposition  gives sufficient conditions for $s$-convex ordering, and is very useful for the verification of higher order convex orders. However, it is worth noticing that in the case of some inequalities, the distribution functions cross more than $s$-times. Therefore a simple application of this proposition is impossible.

In the paper \cite{terRajba2016a}, a theorem on necessary and sufficient conditions for higher order convex stochastic ordering is given. This theorem is a counterpart of the Levin-Ste\v{c}kin theorem \cite{terLevinSteckin1960} concerning convex stochastic ordering. 
Based on this theorem, useful criteria for the verification of higher order convex stochastic ordering  are given. These results can be useful in the study of Hermite-Hadamard type inequalities for higher order convex functions, and in particular inequalities between the quadrature operators.
 It is worth noticing, that these criteria can be easier to checking of higher order convex orders, than those given in \cite{terDenLefSha98, terKARLNOVIK}.

Let  $F_1,F_2 \colon [a,b]\to\mathbb{R}$ be two functions with bounded variation and $\mu_1$, $\mu_2$ be the signed measures corresponding to $F_1$, $F_2 $, respectively. We say that $F_1$ is dominated by $F_2 $ in $(n+1)$-\textit{convex ordering} sense $(n\geq 1)$ if 
\begin{equation*}
\int_{-\infty}^{\infty}f(x)dF_1(x)\leq\int_{-\infty}^{\infty}f(x)dF_2(x)
\end{equation*}
for all $n$-convex functions $f \colon [a,b]\to\mathbb{R}$.
In that case we write $F_1 \leq_{(n+1)-cx} F_2$, or $\mu_1 \leq_{(n+1)-cx} \mu_2$. 
In the following theorem we give necessary and sufficient conditions for $(n+1)$-convex ordering of two functions with bounded variation.
\begin{theorem}[\cite{terRajba2016a}]\label{terth:2.1}
Let $a,b\in\mathbb{R}$, $a<b$, $n\in \mathbb{N}$ and let $F_1,F_2 \colon [a,b]\to\mathbb{R}$ be two functions with bounded variation such that $F_1(a)=F_2(a)$. Then, 
in order that
\begin{equation*}
\int_a^bf(x)dF_1(x)\leq\int_a^bf(x)dF_2(x)
\end{equation*}
for all continuous $n$-convex functions $f \colon [a,b]\to\mathbb{R},$
it is necessary and sufficient that $F_1$ and $F_2$ verify the following  conditions:
\begin{equation*}
F_1(b)= F_2(b),
\end{equation*}
\begin{equation*}
\int_a^b F_1(x)dx = \int_a^b F_2(x)dx, 
\end{equation*}
\begin{multline}\label{tereq:ls2.4} 
\int_a^b \int_a^{x_{k-1}} \ldots \int_a^{x_1} F_1(t)dtdx_1\ldots dx_{k-1}=\\ 
 \int_a^b \int_a^{x_{k-1}} \ldots \int_a^{x_1} F_2(t)dtdx_1\ldots dx_{k-1}\quad\text{for }\: k=2,\ldots ,n,
 \end{multline}
  \begin{multline}\label{tereq:ls2.5} 
(-1)^{n+1}\int_a^x \int_a^{x_{n-1}} \ldots \int_a^{x_1} F_1(t)dtdx_1\ldots dx_{n-1}\leq\\ 
 (-1)^{n+1}\int_a^x \int_a^{x_{n-1}} \ldots \int_a^{x_1} F_2(t)dtdx_1\ldots dx_{n-1} \quad\text{for all }\: x\in(a,b).
 \end{multline}
\end{theorem}
\begin{corollary}[\cite{terRajba2016a}]\label{tercor:2.4}
Let $\mu_1$, $\mu_2$ be two signed measures on $\mathcal{B}(\mathbb{R})$, which are concentrated on $(a,b)$, and such that $\int_a^b\lvert x \rvert ^n \mu_i(dx)<\infty$, $i=1,2$. Then in order that 
\begin{equation*}
\int_a^bf(x)d\mu_1(x)\leq\int_a^bf(x)d\mu_2(x)
\end{equation*}
for continuous $n$-convex functions $f \colon [a,b]\to\mathbb{R}$, it is necessary and sufficient that  $\mu_1$, $\mu_2$ verify the following conditions:
\begin{align}
\mu_1\left((a,b)\right)&=\mu_2\left((a,b)\right),\label{eq:2.17}\\
\int_a^b x  ^k \mu_1(dx)&=\int_a^b x  ^k \mu_2(dx)\quad\text{for }\:\:k=1,\ldots ,n,\label{eq:2.18}\\
\int_a^b \bigl(t-x\bigr)  ^n _+\mu_1(dt)&=\int_a^b \bigl(t-x\bigr) ^n _+\mu_2(dt)\quad\text{for all }\:\:x\in (a,b),\label{eq:2.19}
\end{align}
where $y^n_+=\Bigl\lbrace max \lbrace y,0\rbrace \Bigr\rbrace ^n$, $y\in \mathbb{R}$.
\end{corollary}

In \cite{terDenLefSha98}, can be found the following necessary and sufficient conditions for the verification of the $(s+1)$-convex order.
\begin{proposition}[\cite{terDenLefSha98}]\label{terprop:2.5}
If $X$ and $Y$ are two real valued random variables such that $\mathbb{E} \lvert X \rvert^{s}<\infty$ and  $\mathbb{E} \lvert Y\rvert^{s}<\infty$, then 
\begin{equation*}\label{tereq:2.5c}
\mathbb{E} f(X) \leq \mathbb{E} f(Y)
\end{equation*}
for all continuous $s$-convex functions $f \colon \mathbb{R} \to \mathbb{R}$
if, and only if, 
\begin{align}
\mathbb{E} X^k&=\mathbb{E} Y^k \quad\text{for } \:k=1,2,\ldots,s,\label{tereq:2.5a} \\
\mathbb{E} (X-t)^{s}_+&\leq \mathbb{E} (Y-t)^{s}_+ \quad\text{for all }\: t \in \mathbb{R}\label{tereq:2.5b}.
\end{align}
\end{proposition}
\begin{remark}[\cite{terRajba2016a}]
Note, that if the measures $\mu_X$, $\mu_Y$, corresponding to the random variables $X$, $Y$, respectively, occurring in Proposition \ref{terprop:2.5}, are concentrated on some interval $[a,b]$, then this proposition is an easy consequence of Corollary \ref{tercor:2.4}.
\end{remark}

Theorem \ref{terth:2.1} can be rewritten in the following form.
\begin{theorem}[\cite{terRajba2016a}]\label{terth:2.7a}
Let $F_1,F_2 \colon [a,b]\to\mathbb{R}$ be two functions with bounded variation such that $F_1(a)=F_2(a)$. Let
\begin{align*}
H_0(t_0)&=F_2(t_0)-F_1(t_0) \quad\text{for } \: t_0\in [a,b],\\
H_{k}(t_k)&=\int_a^{t_{k-1}}H_{k-1}(t_{k-1})dt_{k-1}\quad\text{for } \: t_{k}\in [a,b], \  k=1,2,\ldots ,n.
\end{align*}
Then, 
in order that
\begin{equation*}
\int_a^bf(x)dF_1(x)\leq\int_a^bf(x)dF_2(x)
\end{equation*}
for all continuous $n$-convex functions $f \colon [a,b]\to\mathbb{R},$
it is necessary and sufficient that the following  conditions are satisfied:
\begin{align*}
H_k(b)&=0 \quad\text{for } \: k=0,1,2,\ldots ,n,\\
(-1)^{n+1}H_{n}(x)&\geq 0 \quad\text{for all } \: x\in (a,b).
\end{align*}
\end{theorem}
\begin{remark}[\cite{terRajba2016a}]\label{terrmk:2.8a}
The functions $H_1, \ldots,H_n$, that appear in Theorem \ref{terth:2.7a} can be obtained from the following formulas
	\begin{equation}\label{tereq:2.8c}
H_{n}(x)=(-1)^{n+1}\int_a^b \frac{(t-x)_+^n}{n!}d (F_2(t)-F_1(t)), 
\end{equation} 
\begin{equation}\label{tereq:2.8d}
H_{k-1}(x)=H_{k}^{\,'}(x), \quad k=2,3,\ldots ,n.
\end{equation}
\end{remark}

\par\bigskip
Note that the function $(-1)^{n+1}H_{n-1}$, that appears in Theorem \ref{terth:2.7a}, play a role similar to the role of the function $F=F_2-F_1$ in Lemma \ref{terlemma:1.5}. Consequently, from Theorem \ref{terth:2.7a}, Lemma \ref{terlemma:1.5} and Remarks \ref{terrmk:a2}, \ref{terrmk:2.8a}, we obtain immediately the following criterion, which can be useful for the verification of higher order convex ordering.
\begin{corollary}[\cite{terRajba2016a}]\label{tercor:2.7b}
Let $F_1,F_2 \colon [a,b]\to\mathbb{R}$ be functions with bounded variation such that $F_1(a)=F_2(a)$, $F_1(b)=F_2(b)$ and $H_{k}(b)=0$ $(k=1,2,\ldots ,n)$, where $H_k(x)$ $(k=1,2,\ldots ,n)$ are given by \eqref{tereq:2.8c} and \eqref{tereq:2.8d}. Let $a<x_1< \ldots<x_m<b$ be the points of sign changes of the function $H_{n-1}$ and let $(-1)^{n+1}H_{n-1}(x)\geq 0$   for $x\in(a,x_1)$.
\begin{itemize}
\item If $m$ is even then the inequality
 \begin{equation}\label{tereq:2.7c}
\int_a^bf(x)dF_1(x)\leq\int_a^bf(x)dF_2(x),
\end{equation}
is not satisfied by all continuous $n$-convex functions $f \colon [a,b]\to\mathbb{R}$.
 
\item If $m$ is odd, then the inequality \eqref{tereq:2.7c} is satisfied for all continuous $n$-convex functions $f \colon [a,b]\to\mathbb{R}$ if, and only if, 
\begin{equation}\label{tereq:2.7d}
(-1)^{n+1}H_n(x_2)\geq 0,\:\:\; (-1)^{n+1}H_n(x_4)\geq 0,\:\;\ldots,\;\:
 (-1)^{n+1}H_n(x_{m-1})\geq 0.
\end{equation}
\end{itemize}
\end{corollary}

\par\bigskip
In the numerical analysis, some inequalities, which are connected with quadrature operators, are studied. These inequalities, called extremalities, are a particular case of the Hermite-Hadamard type inequalities. Many extremalities are known in the numerical analysis (cf. \cite{terBes08, terBraPet, terBraSch81} and the references therein).
The numerical analysts prove them using the suitable differentiability assumptions. As proved W\k{a}sowicz in the papers \cite{terSzWas07b, terSzWas08, terSzWas10}, for convex functions of higher order, some extremalities can be obtained without assumptions of this kind, using only the higher order convexity itself. The support-type properties play here the crucial role. As we show in \cite{terRajba12f, terRajba12g}, some extremalities can be proved using a probabilistic characterization.The extremalities, which we study are known, however, our method using the Ohlin lemma \cite{terOhlin69} and the Denuit-Lef\`{e}vre-Shaked theorem \cite{terDenLefSha98} on sufficient conditions for the convex stochastic ordering seems to be quite easy. It is worth noticing that, these theorems concern only the sufficient conditions, and they can not be used to the proof some extremalities (see \cite{terRajba12f, terRajba12g}). In these cases, results given in the paper \cite{terRajba2016a}, may be useful .
\par\bigskip

For a function $f:[-1,1]\to\mathbb{R}$ we consider six operators approximating the integral mean
value
\[
 \mathcal{I}(f):=\tfrac{1}{2}\int\limits_{-1}^1f(x)dx.
\]
They are given by
\begin{align*}
 C(f)&:=\tfrac{1}{3}\Bigl(
                      f\bigl(-\tfrac{\sqrt{2}}{2}\bigr)+f(0)+
					  f\bigl(\tfrac{\sqrt{2}}{2}\bigr)
                     \Bigr),\\
 \mathcal{G}_{2}(f)&:=\tfrac{1}{2}
              \Bigl(
			   f\bigl(-\tfrac{\sqrt{3}}{3}\bigr)
               +f\bigl(\tfrac{\sqrt{3}}{3}\bigr)
			  \Bigr),\\
 \mathcal{G}_{3}(f)&:=\tfrac{4}{9}f(0)+
             \tfrac{5}{18}
              \Bigl(
			   f\bigl(-\tfrac{\sqrt{15}}{5}\bigr)+
			   f\bigl(\tfrac{\sqrt{15}}{5}\bigr)
			  \Bigr),\\
 \mathcal{L}_{4}(f)&:=\tfrac{1}{12}\bigl(f(-1)+f(1)\bigr)
              +\tfrac{5}{12}\Bigl(f\bigl(-\tfrac{\sqrt{5}}{5}\bigr)
              +f\bigl(\tfrac{\sqrt{5}}{5}\bigr)\Bigr),\\
 \mathcal{L}_{5}(f)&:=\tfrac{16}{45}f(0)+
              \tfrac{1}{20}\bigl(f(-1)+f(1)\bigr)+
			  \tfrac{49}{180}
			   \Bigl(
			    f\bigl(-\tfrac{\sqrt{21}}{7}\bigr)+
			    f\bigl(\tfrac{\sqrt{21}}{7}\bigr)
			   \Bigr),\\
 S(f)&:=\tfrac{1}{6}\bigl(f(-1)+f(1)\bigr)+\tfrac{2}{3}f(0).
\end{align*}

The operators $\mathcal{G}_{2}$ and $\mathcal{G}_{3}$ are connected with Gauss-Legendre rules. The operators $\mathcal{L}_{4}$ and $\mathcal{L}_{5}$ are connected with Lobatto quadratures. The operators $S$ and $C$ concern Simpson and Chebyshev quadrature rules, respectively. The operator $\mathcal{I}$ stands for the integral mean value (see e.g. \cite{terRal65}, \cite{terWeisCh}, \cite{terWeisLe}, \cite{terWeisLo}, \cite{terWeisSi}).

We will establish all possible inequalities between these operators in the class of higher order convex functions.

\begin{remark}\label{terremark:21}
Let $X_2$, $X_3$, $Y_4$, $Y_5$, $U$, $V$ and $Z$ be random variables such that
\begin{eqnarray*}
\mu_{X_2} & = & \frac{1}{2}\left(\delta_ {- \frac{\sqrt{3}}{3}}+\delta_\frac{\sqrt{3}}{3}\right), \\
\mu_{X_3} & = & \frac{4}{9}\delta_0 + \frac{5}{18}\left(\delta_{-\frac{\sqrt{15}}{5}}+\delta_\frac{\sqrt{15}}{5}\right),\\
\mu_{Y_4} & = & \frac{1}{12}(\delta_{-1}+\delta_1)+\frac{5}{12}\left(\delta_{-\frac{\sqrt{5}}{5}} + \delta_\frac{\sqrt{5}}{5}\right), \\
\mu_{Y_5} & = & \frac{16}{45}\delta_0 + \frac{1}{20}(\delta_{-1}+\delta_1)+\frac{49}{180}\left( \delta_{-\frac{\sqrt{21}}{7}} + \delta_\frac{\sqrt{21}}{7}\right), \\
\mu_U & = & \frac{2}{3}\delta_0 + \frac{1}{6}(\delta_{-1}+\delta_1), \\
\mu_V & = & \frac{1}{3}\left( \delta_{-\frac{\sqrt{2}}{2}} + \delta_0+\delta_\frac{\sqrt{2}}{2}\right), \\
\mu_Z(dx) & = & \frac{1}{2}\chi_{[-1,1]}(x)dx.
\end{eqnarray*}
Then we have
$$
\mathcal{G}_{2}(f) = \mathbb{E}[f(X_2)], \quad \mathcal{G}_{3}(f) = \mathbb{E}[f(X_3)],
$$
$$
\mathcal{L}_{4}(f) = \mathbb{E}[f(Y_4)], \quad \mathcal{L}_{5}(f) = \mathbb{E}[f(Y_5)],
$$
$$
S(f) = \mathbb{E}[f(U)], \quad C(f) = \mathbb{E}[f(V)], \quad \mathcal{I}(f) = \mathbb{E}[f(Z)].
$$
\end{remark}
\begin{theorem}\label{terth:24}
Let $f \colon [-1,1] \to \mathbb{R}$ be 5-convex. Then
\begin{equation}
\mathcal{G}_{3}(f) \leq \mathcal{I}(f)  \leq \mathcal{L}_{4}(f),
\label{tereq:22a}
\end{equation}
\begin{equation}
\mathcal{G}_{3}(f)  \leq \mathcal{L}_{5}(f)\leq \mathcal{L}_{4}(f).
\label{tereq:22}
\end{equation}
\end{theorem}

 Note, that the inequalities \eqref{tereq:22a} and \eqref{tereq:22} can be simply derived from Theorems \ref{terth:20} and \ref{terth:19a} (see \cite{terRajba2016a}).

\begin{remark}\label{terremark:25}
The inequalities \eqref{tereq:22} can be found in \cite{terSzWas08, terSzWas10}. W\k{a}sowicz \cite{terSzWas08} proved, that in the class of 5-convex functions the operators $\mathcal{G}_{2},C,S$ are not comparable both with each other and with $\mathcal{G}_{3},\mathcal{L}_{4},\mathcal{L}_{5}$. 
\end{remark}
\begin{theorem}\label{terth:22}
Let $f \colon [-1,1] \to\mathbb{R}$ be 3-convex. Then
\begin{equation}
\mathcal{G}_{2}(f) \leq \mathcal{I}(f) \leq S(f),
\label{tereq:19}
\end{equation}
\begin{equation}
\mathcal{G}_{2}(f)\leq C(f) \leq T(f) \leq S(f),
\label{tereq:20}
\end{equation}
where $T \in \{ \mathcal{G}_{3}, \mathcal{L}_{5} \}$.
\end{theorem}
\par\bigskip

In \cite{terRajba2016a} is given a new simple proof of Theorem \ref{terth:22}. Note, that from Theorem \ref{terth:20}, we obtain $\mathcal{G}_{3}(f) \leq \mathcal{I}(f)$ and $\mathcal{I}(f) \leq S(f)$, which implies \eqref{tereq:19}. From Theorem \ref{terth:19}, we obtain $\mathcal{G}_{2}(f) \leq C(f)$. By Theorem \ref{terth:19a}, we get $C(f) \leq \mathcal{G}_{3}(f)$, $C(f) \leq \mathcal{L}_{5}(f)$, $\mathcal{G}_{3}(f) \leq S(f)$, $\mathcal{L}_{5}(f)\leq S(f)$. 
\begin{remark}\label{terremark:23}
The inequalities \eqref{tereq:20} can be found in \cite{terSzWas07b}. W\k{a}sowicz \cite{terSzWas07b} proved, that the quadratures $\mathcal{L}_{4}$, $\mathcal{L}_{5}$ and $\mathcal{G}_{3}$ are not comparable in the class of 3-convex functions. 
\end{remark}
\begin{remark}Moreover, W\k{a}sowicz \cite{terSzWas07b, terSzWas08b} proved, that
\begin{equation}
  C(f) \leq \mathcal{L}_{4}(f),
  \label{tereq:20b}
  \end{equation}
if $f$ is 3-convex.

\par\bigskip
 The proof given in \cite{terSzWas07b} is rather complicated. This was done using computer software. In \cite{terSzWas08b}, can be found a new proof of \eqref{tereq:20b}, without the use of any computer software, based on the spline approximation of convex functions of higher order. It is worth noticing, that Proposition \ref{terprop:18} does not apply to proving \eqref{tereq:20b}, because the distribution functions   $F_V$ and $F_{Y_4}$ cross exactly $5$-times.

\par\bigskip
In \cite{terRajba2016a}, the following new proof of \eqref{tereq:20b} is given. In this proof of \eqref{tereq:20b}, we use Corollary \ref{tercor:2.7b}. Note, that we have $F_1=F_V$, $F_2=F_{Y_4}$, $H_0=F=F_{Y_4}-F_V$. By \eqref{tereq:2.8c} and \eqref{tereq:2.8d}, we obtain
\begin{eqnarray*}
H_3(x)&=&\tfrac{1}{72}\left\{ \left( -1-x \right)^3_+ +\left( 1-x \right)^3_+ 
+5\left[\left( -\tfrac{\sqrt{5}}{5}-x \right)^3_+ +\left( \tfrac{\sqrt{5}}{5}-x \right)^3_+\right] \right.\\ \\
\quad &\,& \left.-4\left[\left( -1-x \right)^3_+ +\left( -\tfrac{\sqrt{2}}{2}-x \right)^3_+ +\left( -x \right)^3_+ +\left( \tfrac{\sqrt{2}}{2}-x \right)^3_+ \right] \right\},
\end{eqnarray*}
\begin{eqnarray*}
H_2(x)&=&\tfrac{1}{24}\left\{ -\left( -1-x \right)^2_+ -\left( 1-x \right)^2_+ 
-5\left[\left( -\tfrac{\sqrt{5}}{5}-x \right)^2_+ +\left( \tfrac{\sqrt{5}}{5}-x \right)^2_+\right] \right.\\ \\
\quad &\,& \left.+4\left[\left( -1-x \right)^2_+ +\left( -\tfrac{\sqrt{2}}{2}-x \right)^2_+ +\left( -x \right)^2_+ +\left( \tfrac{\sqrt{2}}{2}-x \right)^2_+ \right] \right\}.
\end{eqnarray*}
Similarly, $H_1(x)$ can be obtained from the equality $H_1(x)=H_2^{\,'}(x)$. We compute, that $x_1=-1-\sqrt{5}+2 \sqrt{2}$, $x_2=0$, $x_3=1+\sqrt{5}-2 \sqrt{2}$ are the points of sign changes of the function $H_2(x)$. It is not difficult to check, that the assumptions of Corollary \ref{tercor:2.7b} are satisfied.  Since
$$  
(-1)^{3+1}H_3(x_2)=(-1)^{3+1}H_3(0)= \tfrac{1}{72}+\tfrac{\sqrt{5}}{360}-\tfrac{\sqrt{2}}{72}>0,
$$
it follows, that the inequalities \eqref{tereq:2.7d} are satisfied. From Corollary \ref{tercor:2.7b} we conclude, that the relation \eqref{tereq:20b} holds.
\end{remark}

\end{document}